\documentclass{article}
\usepackage[utf8]{inputenc}
\usepackage{amsmath}
\usepackage{amsfonts}
\usepackage{amssymb}
\usepackage{amsthm}
\usepackage{soul}
\usepackage{bm}
\usepackage{enumerate}
\usepackage{comment}
\usepackage{graphicx}
\usepackage{float}
\usepackage{hyperref}
\usepackage{arydshln}
\usepackage{mathtools}
\usepackage{bbm}

\usepackage[left=35mm,right=35mm,top=20mm,bottom=20mm,paper=a4paper]{geometry}
\makeatletter
\newcommand*\bigcdot{\mathpalette\bigcdot@{.5}}
\newcommand*\bigcdot@[2]{\mathbin{\vcenter{\hbox{\scalebox{#2}{$\m@th#1\bullet$}}}}}
\makeatother

\newtheorem{Lemma}{Lemma}
\newtheorem{Claim}{Claim}
\newtheorem{Theorem}{Theorem}

\newtheorem{Corollary}{Corollary}

\newcommand{\norm}[1]{\left\lVert#1\right\rVert}

\newsavebox{\overlongequation}
\newenvironment{dontbotheriftheequation*isoverlong}
{\begin{displaymath}\begin{lrbox}{\overlongequation*}$\displaystyle}
 {$\end{lrbox}\makebox[0pt]{\usebox{\overlongequation*}}\end{displaymath}}

\title{The Dirichlet spectrum}
\author{Alon Agin and Barak Weiss }

\begin{document}

\maketitle

\begin{abstract}
\noindent Akhunzhanov and Shatskov defined the Dirichlet spectrum, corresponding to $m \times n$ matrices and to norms on $\mathbb{R}^m$ and $\mathbb{R}^n$. In case $(m,n) = (2,1)$ and using the Euclidean norm on $\mathbb{R}^2$, they showed that the spectrum is an interval. We generalize this result to arbitrary $(m,n) \neq (1,1)$ and arbitrary norms, improving previous works \cite{KR2},\cite{Sch2}\cite{Sch1}\cite{HSW},\cite{AG}. We also define some related spectra and show that they too are intervals. Our argument is a modification of an argument of Khintchine from 1926.    
\end{abstract}
\section{Introduction}

Throughout this paper bold symbols and bold fonts are used for definitions. 

\ \\ Denote by $\boldsymbol{M_{m,n}}=M_{m,n}(\mathbb{R})$
    the space of $m\times n$ matrices with real entries. Let $d=m+n$. \ \\ For $\norm{\cdot}_m$ an arbitrary norm on $\mathbb{R}^m$ and $\norm{\cdot}_n$ an arbitrary norm on $\mathbb{R}^n$, \footnote{\,This notation should not be confused with the $p$-norm defined by $\norm{\vec{x}\,}_p:=\big(\sum|x_i|^p  \big)^{1/p}$.} and for a vector $(\Vec{x},\Vec{y})\in\mathbb{R}^m\times\mathbb{R}^n$ we denote by $\boldsymbol{\norm{\cdot}_{m,n}}$ the norm on $\mathbb{R}^{d}$ defined by

$$\boldsymbol{\norm{(\Vec{x},\Vec{y}\,)}_{m,n}}:=\max\,\{\norm{\Vec{x}\,}_m , \norm{\Vec{y}\,}_n\}.$$ 

\noindent From the Minkowski convex body theorem one deduces a higher dimensional generalization of the
classical Dirichlet’s Diophantine approximation theorem -- there exists a constant $\Delta_0$ such that for every $t\in\mathbb{R}$ sufficiently large and any matrix $\Theta\in M_{m,n}$ there exists $(\Vec{q},\Vec{p}\,)\in\mathbb{Z}^n\times\mathbb{Z}^m$ such that
 \ \\ 
\begin{equation}\label{equation*1}
\Biggl\{
\begin{aligned}
&\,\,0<\norm{\Vec{q}\,}_n\leq t^{1/n}\\
&\,\,t^{1/m}\,\norm{\Theta\Vec{q}-\vec{p}\,}_m<\Delta_0.
\end{aligned} 
\end{equation}

\ \\ Denote by  $\Delta_{\norm{\cdot}_{m,n}}$ the infimum over all $\Delta_0$ which satisfies (\ref{equation*1}). 

\ \\ Let $t_0>0$ be the positive real number which satisfies $t_0=\min\,\{\,\norm{\Vec{q}\,}_n\,|\,\,\Vec{q}\in \mathbb{Z}^n\setminus\{\Vec{0}\}\,\}$. 

\ \\We then define \textbf{the (rescaled) irrationality measure function of \boldsymbol{$M_{m,n}$}} (w.r.t. $\norm{\cdot}_{m,n}$) to be the function $\chi:M_{m,n}\times[t_0,\infty)\longrightarrow[0,\infty)$ defined by

\begin{equation*}
\pmb{\chi(\Theta,t)}=\chi(\Theta,t)_{\norm{\cdot}_{m,n}}:=t^{n/m}\,
\big(\min_{ 0<\norm{\Vec{q}\,}_n\leq \, t\,\,\,\, \Vec{p}\,\in\mathbb{R}^m} \norm{\Theta\Vec{q}-\Vec{p}\,}_m\big).   
\end{equation*}

\ \\ Put in other words, the higher dimensional generalization of Dirichlet's theorem tells us that for every $t$ sufficiently large and every $\Theta\in M_{m,n}$ we have that \ \\  \begin{equation}\label{equation*4}
\chi(\Theta,t)_{\norm{\cdot}_{m,n}}\leq\Delta_{\norm{\cdot}_{m,n}}.
\end{equation}

\ \\ We'll denote until the end of the section the constant $\Delta_{\norm{\cdot}_{m,n}}$ and the function $\chi(\Theta,t)_{\norm{\cdot}_{m,n}}$ only as $\Delta$ and $\chi(\Theta,t)$, keeping in mind these are norm-dependent objects.

\ \\ In case $\limsup_{t\rightarrow\infty} \,\chi(\Theta,t)<\Delta$, we say that $\boldsymbol{\Theta}$ \textbf{is Dirichlet-improvable}, meaning that for this specific $\Theta$, we can improve the upper bound $\Delta$ given to us by Dirichlet's theorem. 

\ \\ A trivial observation is that if $\Theta\in M_{m,n}(\mathbb{Q})$, then for $t$ large enough we have that $\chi(\Theta,t)=0$, so in particular we have that $\limsup_{t\rightarrow\infty} \chi(\Theta,t)=0$. For the typical case, Davenport and Schmidt showed in \cite{DS} that for $n=1$ or $m=1$ and for the maximum norm on $\mathbb{R}^m$ or $\mathbb{R}^n$, respectively, for Lebesgue almost every $\Theta$ we have that $\limsup_{t\rightarrow\infty} \chi(\Theta,t)=\Delta=1$. Extending Davenport and Schmidt's argument (see \cite{KR3}), one can show that for any $m$ and $n$, any arbitrary norms on $\mathbb{R}^m$ and $\mathbb{R}^n$, and for Lebesgue almost every $\Theta\in M_{m,n}$ we have that
\begin{equation}\label{equation*3.333}
\limsup_{t\rightarrow\infty} \chi(\Theta,t)=\Delta.
\end{equation}

\ \\ Following Akhunzhanov and Shatskov in \cite{AS}, define \pmb{the Dirichlet spectrum of $M_{m,n}$} (w.r.t. $\norm{\cdot}_{m,n}$) to be 

\begin{equation*}
\pmb{\mathbb{D}_{m,n}
}=\mathbb{D}_{\norm{\cdot}_{m,n}
}:= \,\big\{\,\limsup\limits_{t\rightarrow\infty} \,\chi(\Theta,t)\,\,|\,\,\Theta\in M_{m,n}
\big\}.
\end{equation*}

\ \\ Again, we denote until the end of the section the Dirichlet spectrum simply by $\mathbb{D}_{m,n}$, keeping in mind this is also a norm-dependent object.

\ \\ By (\ref{equation*4}) we know $\mathbb{D}_{m,n}
\subseteq[0,\Delta]$. In case $n=m=1$, there are many open questions about the spectrum and it appears to have a complicated structure (see e.g \cite{AS}). In particular, we know that $\mathbb{D}_{1,1}\neq[0,\Delta]=[0,1]$. 

\ \\ In this paper we prove the following:

\begin{Theorem}\label{Theorem1}
 Let $m$ and $n$ such that $\max\,(m,n)>1$, let $\norm{\cdot}_m$ and $\norm{\cdot}_n$ be arbitrary norms on $\mathbb{R}^m$ and $\mathbb{R}^n$, and let $\Delta$ be the minimal constant which satisfies (\ref{equation*4}) for every $t$ large enough. Then$$\mathbb{D}_{m,n}=[0,\Delta].$$ 
\noindent Furthermore, for any $c\in[0,\Delta]$, the set $$\{\Theta\in M_{m,n}\,|\,\limsup_{t\rightarrow\infty}\,\chi(\Theta,t)=c\,\}$$ \noindent is uncountable and dense in $M_{m,n}$.
\end{Theorem}

\ \\ Theorem 1 extends many recent results, we present briefly some of them. In \cite{AS}, Akhunzhanov and Shatskov were the first to consider higher dimensions than $m=n=1$. Surprisingly at that time, they showed that for $m=2,n=1$ and for the Euclidean norm on $\mathbb{R}^2$ we have that $\mathbb{D}_{\,2,1}=[0,\Delta]=[0,\sqrt{\frac{2}{\sqrt{3}}}]$. In \cite{Sch1}, Schleischitz showed that for $n=1$, $m\geq2$, and for the maximum norm on $\mathbb{R}^m$, we have that $\mathbb{D}_{m,1}=[0,\Delta]=[0,1]$. In \cite{Sch2}, Schleischitz also showed that for $m=1$, $n\geq2$, and for the maximum norm on $\mathbb{R}^n$, we have that $\mathbb{D}_{1,n}=[0,\Delta]=[0,1]$. In \cite{AG}, the first-named author showed that for $m=2$ and $n=1$ and for any norm on $\mathbb{R}^2$ there exists a positive constant $c$ such that $[0,c]\subseteq \mathbb{D}_{\,2,1} $, with $c=\Delta$ in case $\norm{\cdot}$ is induced from an inner product induced by a diagonal matrix. In \cite{HSW}, Hussain, Schleischitz and Ward showed that for a wide range of integer pairs $(m,n)$ and wide range of norms on $\mathbb{R}^m$ and $\mathbb{R}^n$ there exists a positive constant $c$ such that $[0,c]\subseteq \mathbb{D}_{m,n} $, with $c=\Delta=1$ in case both norms are the maximum norm and if $m|n$ or $n|m$. In \cite{KR2}, Kleinbock and Rao showed in that for $m=2,n=1$ and any norm on $\mathbb{R}^2$, we have that $\Delta$ is an accumulation point of $\mathbb{D}_{\,2,1}$. 

\subsection{Some extensions of the main result}

\ \\ Let $d=m+n$. We say two positive vectors $\Vec{\alpha}\in\mathbb{R}_{>0}^m$ and $\Vec{\beta}\in\mathbb{R}_{>0}^n$ \pmb{are weights on $\mathbb{R}^{d}$} if\, $\sum_{i=1}^{m}\alpha_i=\sum_{i=1}^{n}\beta_i=1$. 

\ \\ Let $\pmb{\mathcal{X}_d}=SL_d(\mathbb{R})/SL_d(\mathbb{Z})$ denote the space of unimodular lattices in $\mathbb{R}^d$.   Given $
\Theta\in M_{m,n}$ we denote by $\Lambda_\Theta$ the unimodular lattice defined by

$$\pmb{\Lambda_{\Theta}}:=\begin{pmatrix}
I_m & \Theta \\
0 & I_n \\
\end{pmatrix}\mathbb{Z}^d.$$

\ \\ Given $t>0$ and $\vec{\alpha},\vec{\beta}$ weights on $\mathbb{R}^{d}$ we define
$$\pmb{g_{t,\vec{\alpha},\vec{\beta}\,}}:=diag\,(t^{\alpha_1},\,...\,,t^{\alpha_m},t^{-\beta_1},\,...\,,t^{-\beta_n})\in SL_d(\mathbb{R}).$$

\noindent Let $b\in(0,\infty]$ and let $\{E_\varepsilon\,|\,\varepsilon\in (0,b]\,\}$ be a continuous decreasing exhaustion of $\mathcal{X}_d$ by compact sets (see (\ref{list}), section 3). 

\ \\Motivated by a question posed by Dmitry Kleinbock, for $\Lambda\in\mathcal{X}_d$ we define \pmb{the extended } \pmb{Dirichlet constant of $\Lambda$} (w.r.t. $E_\varepsilon,\vec{\alpha},\vec{\beta}$\,) to be 

$$ \pmb{dir(\Lambda)\,(E_{\varepsilon},\vec{\alpha},\vec{\beta})}:=\inf\,\{ \,\varepsilon>0\,\,|\,\,g_{t,\vec{\alpha},\vec{\beta}\,} \Lambda\notin\,E_\varepsilon\,\,for\,\,all\,\, t\,\,large\,\,enough\}, $$ with the convention that $\inf\varnothing=b$. Now define \pmb{the extended Dirichlet spectrum of $\mathcal{X}_d$} (w.r.t. $E_\varepsilon,\vec{\alpha},\vec{\beta}$\,) to be 
$$
\pmb{\mathbb{D}_{d}\,(E_\varepsilon,\vec{\alpha},\vec{\beta})}:=\{\, dir(\Lambda)\,(E_{\varepsilon},\vec{\alpha},\vec{\beta}) \,|\, \Lambda\in X_d\,\}.
$$ 

\noindent In this paper we also prove the following:

\begin{Theorem}
Let $d\geq3$. Let $\{E_\varepsilon\,|\,\varepsilon\in (0,b]\,\}$ be a continuous decreasing exhaustion by compact sets of $\mathcal{X}_d$ for some $b\in(0,\infty]$. Let $m,n\in\mathbb{N}$ with $m+n=d$ and let $\Vec{\alpha}$ and $\Vec{\beta}$ be arbitrary weights on $\mathbb{R}^{d}$. Then $$\mathbb{D}_{d}\,(E_\varepsilon,\vec{\alpha},\vec{\beta})=[0,b].$$ \noindent Furthermore, for any $c\in [0,b]$ and any open set $V\in M_{m,n}$ we have that the set  
$$\{\Theta\in V\,|\,dir(\Lambda_\Theta)\,(E_{\varepsilon},\vec{\alpha},\vec{\beta})=c\,\}$$ \noindent is uncountable. 
\end{Theorem}

\noindent For a more detailed discussion see section 3, during which we also explain why Theorem 2 is a generalization of Theorem 1.

\ \\  The paper is structured as follows. In the following section we generalise the setting by changing the function $t^{-1}$ to a more general one, still staying within the realm of Diophantine approximation. For $\Theta\in M_{m,n}$, we discuss $\hat{\omega}(\Theta)$ the uniform exponent of $\Theta$, and use a generalisation of Theorem 1 to prove existence of matrices exhibiting special properties with respect to their uniform exponent (Theorem 3 and Corollary 1). In section 3, we generalise the setting to the space of unimodular lattices in more detail, turning the question to one which focuses on the behaviour of certain trajectories under one-parameter diagonalizable flows (Theorem 2 and 4). In sections 4 and 5 we show that Theorems 2 and 3 follow from a more general topological theorem (Theorem 5). Theorem 5 is the main result of this paper. Throughout the paper we use some technical lemmas (Lemmas 1-4); 
For readability, we postpone the proofs of these lemmas to a separate and final section.

\ \\ The origin of the topological theorem (Theorem 5) goes back to a famous work by Khintchine from the 20's \cite{KHI} regarding existence of non-obvious singular $2$-dimensional vectors, which was later generalised and modified by many authors (see \cite{CAS},\cite{DANI},\cite{WEISS}, \cite{KMW}). This new version of it is a modification of a modification introduced by the second-named author in \cite{WEISS} (Theorem 2.5).

\ \\ \textbf{Acknowledgements.} This work is part of the first author’s doctoral dissertation, conducted at Tel Aviv University under the supervision of the second author. We are grateful to Dmitry Kleinbock and Anurag Rao for useful discussions. The support of grants ISF-NSFC 3739/21 and ISF 2021/24 are gratefully acknowledged.

\section{$\psi$-Dirichlet spectrum and uniform exponent}

\noindent For $\gamma>0$ we define \pmb{the $\gamma$-irrationality measure function of $M_{m,n}$} (w.r.t. $\norm{\cdot}_{m,n}$) to be the function $\chi_\gamma:M_{m,n}\times[t_0,\infty)\longrightarrow[0,\infty)$ defined by

\begin{equation*}
\pmb{\chi_{\gamma}(\Theta,t)}=\chi_\gamma(\Theta,t)_{\norm{\cdot}_{m,n}}:=t^\gamma\,
\big(\min_{ 0<\norm{\Vec{q}\,}_n\leq \, t\,\,\,\, \Vec{p}\,\in\mathbb{R}^m} \norm{\Theta\Vec{q}-\Vec{p}\,}_m\big).
\end{equation*}

\noindent One can easily show that if $\limsup_{t\rightarrow\infty}\chi_{\gamma}(\Theta,t)=c$ for some $c\in(0,\infty)$ then for any $\delta,\varepsilon$ with $0<\delta<\gamma<\varepsilon$ we have $\limsup_{t\rightarrow\infty}\chi_\delta(\Theta,t)=0$ and $\limsup_{t\rightarrow\infty}\chi_{\varepsilon}(\Theta,t)=\infty$. So by Dirichlet's theorem we have that if $\gamma<n/m$ then $\limsup_{t\rightarrow\infty}\chi_{\gamma}(\Theta,t)=0$ for every $\Theta\in M_{m,n}$.

\ \\ Now for $\Theta\in M_{m,n}$ define \pmb{the uniform exponent of $\Theta\in M_{m,n}$} (w.r.t. $\norm{\cdot}_{m,n}$) to be

$$\pmb{\hat{\omega}(\Theta)}=\hat{\omega}(\Theta)_{\norm{\cdot}_{m,n}}:=\sup\big\{\,\gamma>0\,|\,\limsup\limits_{t\rightarrow\infty} \,\chi_{\gamma}(\Theta,t)<\infty\,\big\}.$$

\ \\ By Dirichlet's theorem, clearly we have $\hat{\omega}(\Theta)\geq n/m$ for every $\Theta\in M_{m,n}$, and by (\ref{equation*3.333}) we have $\hat{\omega}(\Theta)= n/m$ for Lebesgue almost every $\Theta\in M_{m,n}$  (see e.g \cite{KMW} for further discussion and results).

\ \\ More generally, let $t>0$ and let $\psi(t)$ be a positive continuous decreasing function. \\ For $\Theta\in M_{m,n}$ fixed define

\begin{equation}\label{equation*5.1111}
\pmb{\lambda_{\Theta,\psi}(t)}=\lambda_{\Theta,\psi}(t)_{\norm{\cdot}_{m,n}} 
:=\min_{(\vec{q},
\vec{p}\,)\in\,\mathbb{Z}^n\times\mathbb{Z}^m\setminus\{\vec{0}\}}\,\max\, \Biggl\{
\begin{array}{l}
t^{-1/n}\,\norm{\Vec{q}\,}_n
   \vspace*{0.2cm}\\
\psi(t)^{-1/m}\,\norm{\Theta\Vec{q}-\vec{p}\,}_m\,
\end{array}
\end{equation}

\ \\ and for $\gamma>0$ define $\pmb{\psi_\gamma(t)}:=t^{-\gamma}$. 

\ \\ That is \begin{equation*}
\lambda_{\Theta,\psi_\gamma}(t)=\min_{(\vec{q},
\vec{p}\,)\in\,\mathbb{Z}^n\times\mathbb{Z}^m\setminus\{\vec{0}\}}\,\max\, \Biggl\{
\begin{array}{l}
t^{-1/n}\,\norm{\Vec{q}\,}_n
   \vspace*{0.2cm}\\
t^{\,\gamma/m}\,\norm{\Theta\Vec{q}-\vec{p}\,}_m\,
\end{array}.
\end{equation*}

\ \\For $\Theta$ fixed, clearly $\lambda_{\Theta,\psi_\gamma}(t)$ and $\chi_\gamma(\Theta,t)$ are closely related functions; note however that $\lambda_{\Theta,\psi_\gamma}(t)$ is a continuous function of $t$, but $\chi_\gamma(\Theta,t)$ is not. In section 6 we give some elementary statements connecting these two functions. In particular, we prove Lemma \ref{Lemma6.111} which states that for $\gamma>0$ we have that
\begin{equation}\label{equation*6.11111}
\limsup_{t\rightarrow\infty}\,\lambda_{\Theta,\psi_\gamma}(t)\,^{1+\gamma n/m}
=\limsup_{t\rightarrow\infty} \,\chi_{\gamma n/m}(\Theta,t).
\end{equation}

\ \\ Now for $\psi(t)$ positive continuous and decreasing, we define \pmb{the $\psi$-Dirichlet spectrum} \\ \pmb{of $M_{m,n}$} (w.r.t. $\norm{\cdot}_{m,n}$) to be 
$$
\pmb{\mathbb{D}_{m,n}(\psi)}:=\big\{\limsup\limits_{t\rightarrow\infty} \,\lambda_{\Theta,\psi}(t)\,\,|\,\,\Theta\in M_{m,n}
\big\}.$$

\noindent In particular  for $\gamma>0$ we have that \begin{equation*}
\mathbb{D}_{m,n}(\psi_{\gamma m/n})=\big\{\limsup_{t\rightarrow\infty} \,\chi_{\gamma}(\Theta,t)\,\,|\,\,\Theta\in M_{m,n}
\big\}. 
\end{equation*}

\ \\ In this paper we also prove the following:

\begin{Theorem}\label{Theorem2}
Let $m$ and $n$ such that $\max\,(m,n)>1$, let $\norm{\cdot}_m$ and $\norm{\cdot}_n$ be arbitrary norms on $\mathbb{R}^m$ and $\mathbb{R}^n$, and let $\Delta$ be the minimal constant which satisfies (\ref{equation*4}). Let $\psi(t)$ be a positive continuous decreasing function such that $\psi(t)=o(t^{-1})$. In case $n=1$, assume additionally that $\lim_{t\rightarrow\infty}t^{\,m}\,\psi(t)=\infty$. Then $$\mathbb{D}_{m,n}\,(\psi)=[0,\infty].$$
\noindent Furthermore, for any $c\in[0,\infty]$, the set $$\{\Theta\in M_{m,n}\,|\,\limsup_{t\rightarrow\infty} \,\lambda_{\Theta,\psi}(t)=c\,\}$$\noindent is uncountable and dense in $M_{m,n}$.
\end{Theorem}

\newpage \noindent An application of Theorem 3 is the following:

\begin{Corollary}
Let $m$ and $n$ such that $\max\,(m,n)>1$, and let $\norm{\cdot}_m$ and $\norm{\cdot}_n$ be arbitrary norms on $\mathbb{R}^m$ and $\mathbb{R}^n$. Let $\gamma>n/m$. In case $n=1$ assume additionally $\gamma<1$. 

\ \\Then for every $c\in[0,\infty]$ there exists uncountably many $\Theta\in M_{m,n}$ such that $\hat{\omega}(\Theta)=\gamma$ and $\limsup_{t\rightarrow\infty} \,\chi_{\gamma}(\Theta,t)=c$. \,\,Furthermore, the set $$\{\Theta\in M_{m,n}\,|\,\hat{\omega}(\Theta)=\gamma\,\,and\,\,\limsup_{t\rightarrow\infty} \,\chi_{\gamma}(\Theta,t)=c\,\}$$ \noindent is uncountable and dense in $M_{m,n}$.
\end{Corollary}

\noindent Corollary 1 follows immediately from Theorem 3 and formula (\ref{equation*6.11111}), and the proof is left for the reader.  

\section{One-parameter diagonalizable flows on $\mathcal{X}_d$}

 \ \\ Let $\mathcal{X}_d$,
$\Lambda_\Theta$, $\vec{\alpha},\vec{\beta}$ and $g_{t,\vec{\alpha},\vec{\beta}\,}$ as in section 1.1. Let $\norm{\cdot}_d$ be an arbitrary norm on $\mathbb{R}^d$.

\ \\ For convenience, denote $\pmb{\Vec{m}_1}:=(1/m\,,\,...\,,\,1/m)\in\mathbb{R}^m$ and $\pmb{\Vec{n}_1}:=(1/n\,,\,...\,,\,1/n)\in\mathbb{R}^n$.

\ \\ For a lattice $\Lambda$ (not necessarily unimodular) define \pmb{the Minkowski's first} \pmb{successive minimum} (w.r.t. $\norm{\cdot}_d$) to be $$
\pmb{\lambda_{1}(\Lambda)}=\lambda_{1}(\Lambda)(\norm{\cdot}_d):=\inf\{\,\norm{\vec{u}}_d\,|\,\vec{0}\neq\vec{u}\in \Lambda\}=\min\{\,\norm{\vec{u}}_d\,|\,\vec{0}\neq\vec{u}\in \Lambda\}    
$$ where the last equality holds by discreteness. Most of the time we'll denote the function $\lambda_1(\norm{\cdot}_d)$ only as $\lambda_1$, keeping in mind that this is a norm-dependent function.

\ \\ If $\norm{\cdot}_d=\norm{\cdot}_{m,n}$ for two arbitrary norms on $\mathbb{R}^m$ and $\mathbb{R}^n$, and if $\vec{\alpha}=\vec{m}_1$ and $\vec{\beta}=\vec{n}_1$, then for any $\Theta\in M_{m,n}$ and any $t>0$ a simple calculation shows that  

\begin{equation}\label{equation*10.111}
\lambda_1(g_{t,\vec{m}_1,\vec{n}_1\,}\Lambda_\Theta)=\lambda_{\Theta,t^{-1}}(t)= \min_{(\vec{q},
\vec{p}\,)\in\,\mathbb{Z}^n\times\mathbb{Z}^m\setminus\{\vec{0}\}}\,\max\, \Biggl\{
\begin{array}{l}
t^{-1/n}\,\norm{\Vec{q}\,}_n
   \vspace*{0.2cm}\\
t^{1/m}\,\norm{\Theta\Vec{q}-\vec{p}\,}_m\,\,
\end{array}
\end{equation}

\noindent for $\lambda_{\Theta,t^{-1}}$ as in formula (\ref{equation*5.1111}). So by Lemma \ref{Lemma6.111} (with $\gamma=1$) we have that

\begin{equation}\label{equation*6.1111}
\limsup\limits_{t\rightarrow\infty}\,\lambda_1(g_{t,\vec{m}_1,\vec{n}_1\,}\Lambda_\Theta)\,^{1+n/m}=\limsup\limits_{t\rightarrow\infty} \,\chi(\Theta,t)_{\norm{\cdot}_{m,n}}. 
\end{equation} \noindent 

\ \\ So with a little more work, one can show that in the specific case of $\norm{\cdot}_d=\norm{\cdot}_{m,n}$, for any $\Theta\in M_{m,n}$ and any $t>0$ we have that 
\begin{equation*}
\lambda_1(g_{t,\vec{m}_1,\vec{n}_1\,}\Lambda_\Theta)\,^{1+n/m}\leq \Delta_{\norm{\cdot}_{m,n}},
\end{equation*} 
\noindent for $\Delta_{\norm{\cdot}_{m,n}}$ as defined after (\ref{equation*1}). 

\ \\ Furthermore, one can show that in fact for any $\Lambda\in \mathcal{X}_d$ we still have that

\begin{equation*}
\lambda_1(\Lambda)\leq \Delta_{\norm{\cdot}_{m,n}}^{m/(m+n)}
\end{equation*} \noindent and that $\Delta_{\norm{\cdot}_{m,n}}$ is the minimal constant which satisfies this inequality.

\ \\ Notice that for general norm and weights $\norm{\cdot}_d,\Vec{\alpha},\Vec{\beta}$, the value of the function $\lambda_1(g_{t,\vec{\alpha},\vec{\beta}\,}\Lambda_\Theta)$ cannot be translated into a set of formulas similar to (\ref{equation*10.111}).

\ \\ Nevertheless, we still have that there exists a minimal constant $r_{\,\norm{\cdot}_d}$ such that for any $\Lambda\in \mathcal{X}_d$ we have that 
\begin{equation}\label{equation*9.1111}
\lambda_1(\Lambda)\leq r_{\,\norm{\cdot}_d},\end{equation}

\noindent and in case $\norm{\cdot}_d=\norm{\cdot}_{m,n}$ we have that 

\begin{equation*}\label{equation*10.1111}
r_{\,\norm{\cdot}_{m,n}}=\Delta_{\norm{\cdot}_{m,n}}^{m/(m+n)}. 
\end{equation*}

\ \\ Now for $\varepsilon>0$ define 
$\pmb{K_\varepsilon}=K_{\varepsilon}(\norm{\cdot}_d):=\{\Lambda\in \mathcal{X}_d\,|\,\varepsilon\leq\lambda_{1}(\Lambda)\,\}$. As before, we'll denote this set only as $K_\varepsilon$, keeping in mind that this is also a norm-dependent set. By Mahler's compactness criterion and by (\ref{equation*9.1111}) we have that $\{\,K_\varepsilon\,|\,0<\varepsilon\leq r_{\,\norm{\cdot}_d} \,\}$ is a continuous decreasing exhaustion by compact sets of $\mathcal{X}_d$.

\ \\With this notation, for $\Lambda\in \mathcal{X}_d$ we have that 
\begin{equation}\label{equation*10.1111}
\begin{split}
\lambda_1(\Lambda)&=\inf\,\,\{\varepsilon>0\,\,|\,\,\Lambda \notin K_\varepsilon\,\},
\end{split}
\end{equation}
\noindent  and in particular for $\Theta\in M_{m,n}$ we have that 
\begin{equation}\label{equation*16.33}
\begin{split}
\limsup_{t\rightarrow\infty}\lambda_1(g_{t,\vec{\alpha},\vec{\beta}\,}\,\Lambda_{\Theta})=\inf\,\,\{ \,\varepsilon>0\,\,|\,\,g_{t,\vec{\alpha},\vec{\beta}\,}\,\Lambda_{\Theta}\notin K_\varepsilon\,for\,\,all\,\, t\,\,large\,\,enough\},
\end{split}
\end{equation}

\noindent with the convention that $\inf\varnothing=r_{\norm{\cdot}_d}$.

\ \\ So in other words, Theorem 1 tells us that for every $c\in[0,r_{\,\norm{\cdot}_{m,n}}]$ there exists uncountably many $\Theta\in M_{m,n}$, corresponding to unimodular lattices of the form $\Lambda_{\Theta}$, such that \begin{equation*}\label{equation*17.111} 
c=\inf\,\{ \,\varepsilon>0\,\,|\,\,g_{t,\vec{m}_1,\vec{n}_1\,}\,\Lambda_{\Theta}\notin K_\varepsilon\,for\,\,all\,\, t\,\,large\,\,enough\}
\end{equation*}

\noindent for any arbitrary norms and weights.

\ \\ More generally, given $b$ with $0<b\leq\infty$ and a collection of subsets of $\mathcal{X}_d$ $\{E_\varepsilon\,|\,\varepsilon\in (0,b]\,\}$, we say that \pmb{$\{E_\varepsilon\,|\,\varepsilon\in (0,b]\,\}$ is a continuous decreasing exhaustion of $\mathcal{X}_d$ by} \pmb{compact sets} if 
 \begin{itemize}

\item $E_\varepsilon$ is compact for any $\varepsilon\in(0,b)$.
\item If $\varepsilon_1<\varepsilon_2$ then $E_{\varepsilon_2}\subset E_{\varepsilon_1}$.
\item  $int\,E_\varepsilon\neq\varnothing$ for any $\varepsilon<b$, $int \,E_b=\varnothing$.
\item Any compact subset of $\mathcal{X}_d$ is contained in $E_\varepsilon$ for some $\varepsilon>0$.
    \item The map $\Lambda\mapsto\inf\,\,\{\varepsilon>0\,\,|\,\, \Lambda\notin\,E_\varepsilon\,\}$ is continuous (where in case $b=\infty$ we endow $[0,\infty]$ with the one point compactification topology). 
\end{itemize}\begin{equation}\label{list}
\end{equation}

\noindent For simplicity we will focus on decreasing exhaustions, though an analogous increasing case is identical up to some obvious modifications. For example, given $\rho(\cdot,\cdot)$ a metric on $\mathcal{X}_d$ and a lattice $\Lambda_0$, we denote by $\pmb{B_{\rho}(\Lambda_0,\varepsilon)}$ the open ball w.r.t. $\rho$ around $\Lambda_0$ of radius $\varepsilon$. So in particular if $\rho$ is a  proper metric, then $\{\,\overline{B_{\rho}(\Lambda_0,\varepsilon)}\,|\,0\leq\varepsilon<\infty \,\}$ is an increasing and continuous exhaustion by compact sets of $\mathcal{X}_d$ for any $\Lambda_0\in \mathcal{X}_d$.

\ \\ By (\ref{equation*6.1111}) and (\ref{equation*10.1111}), Theorem 1 immediately follows from Theorem 2 as the special case $E_\varepsilon=K_\varepsilon(\norm{\cdot}_{m,n}),\vec{\alpha}=\vec{m}_1,\vec{\beta}=\vec{n}_1$. Furthermore, Theorem 2 also allows us to easily deduce the following: 
\begin{Corollary}
Let $d\geq3$. Let $\rho(\cdot,\cdot)$ be a proper metric on $\mathcal{X}_d$, and let $\Lambda_0\in \mathcal{X}_d$ be arbitrary. Let $m,n\in\mathbb{N}$ with $m+n=d$ and let $\Vec{\alpha}$ and $\Vec{\beta}$ be arbitrary weights on $\mathbb{R}^{d}$. 

\ \\ Then for every $c\in[0,\infty]$, the set $$\{\Theta\in M_{m,n}\,|\,\sup\,\{ \,\varepsilon\geq0\,\,|\,\,g_{t,\vec{\alpha},\vec{\beta}\,}\,\Lambda_{\Theta}\notin 
 B_{\rho}(\Lambda_0,\varepsilon)\,for\,\,all\,\, t\,\,large\,\,enough\}=c\,\}$$ \noindent is uncountable and dense in $M_{m,n}$ (with the convention that $\sup\varnothing=0$).
\end{Corollary}

\begin{Corollary}
Let $d\geq3$. Let $\norm{\cdot}_d$ be an arbitrary norm on $\mathbb{R}^d$, and let \,$r_{\norm{\cdot}_d}$ be the minimal constant which satisfies (\ref{equation*9.1111}) for every $\Lambda\in\mathcal{X}_d$. Let $m,n\in\mathbb{N}$ with $m+n=d$ and let $\Vec{\alpha}$ and $\Vec{\beta}$ be arbitrary weights on $\mathbb{R}^{d}$. 

\ \\ Then for every $c\in[0,r_{\norm{\cdot}_d}]$, the set $$\{\Theta\in M_{m,n}\,|\,\limsup_{t\rightarrow\infty}\lambda_1(g_{t,\vec{\alpha},\vec{\beta}\,}\,\Lambda_{\Theta})=c\,\}$$ \noindent is uncountable and dense in $M_{m,n}$.
\end{Corollary}

\noindent Corollary 2 answers a question of Kleinbock (unpublished). Corollary 3 answers a question of Hussain, Schleischitz and Ward posed in \cite{HSW}.  

\ \\ One can now ask whether we have an analogue of Theorem 2 or Corollary 2 in case $d=2$, noticing that although $\mathbb{D}_{1,1}\neq[0,1]$, in the extended Dirichlet spectrum of $\mathcal{X}_2$ we also vary arbitrary lattices which are not of the form $\Lambda_{\Theta}$ for $\Theta\in M_{1,1}$.

\ \\ In the language of Corollary 2 for example, given  $\Lambda_0\in \mathcal{X}_2$ and $\rho(\cdot,\cdot)$ a proper metric on $\mathcal{X}_2$, one can still ask whether for each $c\in[0,\infty]$ there exists $\Lambda\in \mathcal{X}_2$ such that 
\begin{equation*}\label{equation*29.11111}
c=\sup\,\{\,\varepsilon\geq0\,|\,diag(t,t^{-1})\Lambda\notin B_{\rho}(\Lambda_0,\varepsilon)\,\,for\,\,all\,\,t\,\,large\,\,enough\,\}, \end{equation*}

\noindent noticing that in this case there is only one choice for the weights on $\mathbb{R}^2$ and instead of $\mathbb{D}_{2}\,(E_\epsilon, 1,1)$ we will simplify the notation and write $\mathbb{D}_2(E_\epsilon)$.

\ \\ We answer this question in the negative by proving the following:

\begin{Theorem}
Let $d=2$. Let $\{E_\varepsilon\,|\,\varepsilon\in (0,b]\,\}$ be a continuous decreasing exhaustion by compact sets of $\mathcal{X}_2$ for some $b\in(0,\infty]$. Then $$\mathbb{D}_{2}\,(E_\varepsilon)\neq[0,b].$$ \noindent Furthermore, there exists $r\in(0,b)$ such that $(0,r)\cap \mathbb{D}_{2}\,(E_\varepsilon)=\varnothing$.
\end{Theorem}
\noindent \textbf{\textit{Proof of Theorem 4.}} \,\,\,Indeed, a known fact (see e.g. proof of Lemma 11.29 in Einsiedler-Ward, Ergodic Theory with a view toward Number Theory) is that in $\mathcal{X}_2$ there exists a compact $K$ such that any nondivergent orbit intersects $K$ infinitely often. Take $r>0$ such that $K\subset E_r$.
So if the trajectory $diag(t,t^{-1})\Lambda$ is divergent then $dir(\Lambda)(E_{\varepsilon})=0$, and if the trajectory is nondivergent then $dir(\Lambda)(E_{\varepsilon})$ is at least $r$. Thus $\mathbb{D}_{2}\,(E_\varepsilon)$ does not contain any $c\neq 0$ with $c<r$. 
\qed

\section{A topological theorem}

\noindent Given a topological space $\Xi$ and $(X_n)_{n=1}^{\infty}$ a sequence of subsets of $\Xi$, and given $\Omega$ an open set in $\Xi$, we say \pmb{$\Omega$ is connected via $(X_n)$} if there exists $j\in\mathbb{N}$ such that $\Omega\cap X_j\neq\O$ and such that for any $y\in \Omega$ and any $U_y$ a neighborhood of $y$, there exists $i\in\mathbb{N}$ such that $\Omega\,\cap\, X_i$ is connected, $\Omega\cap X_i\cap U_y\neq\O$ and $\Omega\cap X_i\cap X_j\neq\O$.

\ \\ We say \pmb{the space $\Xi$ is locally connected via $(X_n)$} if for any $x\in \Xi$ and every $U_x$ a neighborhood of $x$ there exists $\Omega_x$ an open neighborhood of $x$ with compact closure such that $\Omega_x$ is connected via $(X_n)$ and satisfies $\overline{\Omega_x}\subset U_x$. 

\ \\ Notice that in particular if $\Xi$ is locally connected via $(X_n)$ then $\Xi$ is locally compact, and that
$\cup_{n}X_n$ is dense in $\Xi$.

\ \\ We now prove the following: 
\begin{Theorem}
Let $\Xi$ be a Hausdorff space.

\ \\ Let $f:\Xi\times(0,\infty)\longrightarrow[0,\infty]$ be a continuous function, where we endow $[0,\infty]$ with the one point compactification topology.

\ \\ Let 
$(X_n)_{n=1}^{\infty}$ be a sequence of subsets of $\Xi$ such that the following hold:

\begin{itemize}
\item The space $\Xi$ is locally connected via $(X_n)$.
\item \pmb{Upper uniformity of $(X_n)$ w.r.t. $f$}: For every $i\in\mathbb{N}$ and every $x\in X_i$, $\lim_{t\rightarrow\infty}f(x,t)=0$ uniformly on compact sets. Namely, for every $c>0$, every $i\in\mathbb{N}$, and every compact set $K\subseteq\Xi$ there exists $T\in\mathbb{R}$ such that if $t\geq T$ then $f(x,t)<c$ for all $x\in X_i\cap K$.
\end{itemize} 

\noindent Assume that there exists $0<b\leq\infty$ such that the set $\{x\in\Xi\,|\,\limsup_{t\rightarrow\infty}$ $f(x,t)\geq b$\,\} is dense in $\Xi$.

\ \\ Then for every non-empty open set $V\in\Xi$  and for all $c$ with $0<c<b$ there exists uncountably many $x_0\in V$ such that $\limsup_{t\rightarrow\infty}f(x_0,t)=c$, and in addition $f(x_0,t)<c$ for all large enough $t$.
\end{Theorem}

\ \\ \noindent \textbf{\textit{Proof of Theorem 5.}} \,\,\, Let $c\in(0,b)$, and let 
$V$ be an open non-empty set in $\Xi$. Let $(\varepsilon_k)_{k=0}^{\infty}$ be any decreasing sequence of positive real numbers with $\lim_{k\rightarrow\infty}\varepsilon_k=0$. \ \\ We construct inductively a sequence of non-empty open sets $(\Omega_k)_{k=0}^{\infty}$ with compact closure in $\Xi$, two unbounded increasing sequences of real numbers $(t_k)_{k=0}^{\infty}$, $(s_k)_{k=1}^{\infty}$, and a sequence of natural indices $(j_k)_{k=0}^{\infty}
$ which satisfy the following:

\begin{enumerate}[(A)]
\item $\Omega_k$ is connected via $(X_n)$ for all $k\geq0$.
\item $\overline{\Omega_0}\subset V$ and $\overline{\Omega_{k}}\subset\Omega_{k-1}$ for all $k\geq1$.
\item For all $k\geq1$, all $x\in\Omega_k$ and for all $t$ with $t_{k-1}\leq t\leq t_k$ we have $f(x,t)<c$.    
\item For all $k\geq1$ we have $s_k\in[t_{k-1},t_k]$, and for all $x\in\Omega_k$ we have $c-\varepsilon_k<f(x,s_k)$. 
\item For all $k\geq0$ we have $X_{j_k}\cap\Omega_k\neq\varnothing$, and for all $x\in X_{j_k}\cap\Omega_k$ and all $t\geq t_k$ we have $f(x,t)<c-\varepsilon_k$.
\end{enumerate}

\noindent First we show why sequences which satisfy (B),(C),(D) suffice to get a point $x_0$ which satisfies $\limsup_{t\rightarrow\infty}f(x_0,t)=c$, noticing that conditions (A),(E) will be needed only for the construction of the next step in the inductive procedure. 

\ \\ $\Omega_k$ is non-empty and open with compact closure for all $k$, hence the intersection $\cap_k \Omega_k$ is not empty by condition (B). Let $x_0\in\cap_k \Omega_k$. By condition (C), $\limsup_{t\rightarrow\infty}f(x_0,t)\leq c$, and $f(x_0,t)<c$ for all $t>t_0$. Since $\lim_{k\rightarrow\infty}\varepsilon_k=0$, we get by condition (D) that $\limsup_{t\rightarrow\infty}f(x_0,t)=c$.

\ \\ Now let us construct the sequences inductively. By locally connectedness of $\Xi$ via $(X_n)$, choose $\Omega_0$ open, non-empty with compact closure and connected via $(X_n)$ which satisfies $\overline{\Omega_0}\subset V$, and choose $j_0\in\mathbb{N}$ to be any $j$ such that $X_{j_0}\cap\Omega_0\neq\varnothing$. By upper uniformity, let $T_{0}>0$ such that if $t\geq T_{0}$ then $f(x,t)<c-\varepsilon_0$ for all $x\in X_{j_0}$. Set $t_0=T_{0}$. So $
\overline{\Omega_0}\subset V$, conditions (A),(E) are satisfied for the base of induction, and assuming both are satisfied for $1,2,...,k$ we show how to continue the construction for the $k+1$ step (we do not need to assume inductively properties (B),(C),(D) hold for the $k$th step in order to do so).

\ \\ By the assumption that the set $\{x\in\Xi\,|\,\limsup_{t\rightarrow\infty}$ $f(x,t)\geq b$\,\} is dense in $\Xi$ and since $c<b$, there exist $y\in\Omega_k$ and $s\geq t_k$ arbitrarily large such that $c<f(y,s)$. By continuity of $f$ in $\Xi$, choose $U_y\subset \Omega_k$ to be a small enough neighborhood of $y$ such that for all $x\in U_y$ we still have that $c<f(x,s)$.

\ \\ Let $j_k\in\mathbb{N}$ which satisfies (E) for the kth step. By the inductive assumption $\Omega_k$ is connected via $(X_n)$, so there exists an index $i\in\mathbb{N}$ such that $\Omega_k\cap X_i$ is connected, $\Omega_k\cap X_i\cap U_y\neq\varnothing$ and $\Omega_k\cap X_i\cap X_{j_k}\neq\varnothing$.

\ \\ $\Omega_k$ has compact closure, so by upper uniformity, let $T\in\mathbb{R}$ such that $f(x,t)<c-\varepsilon_{k+1}$ for all $t\geq T$ and all $x\in \Omega_k\cap X_i$.

\ \\ Define $t_{k+1}=\max\,\{s+1,T\}$, and let $\widetilde{y}\in \Omega_k\cap X_i\cap U_y$ and $\widetilde{z}\in \Omega_k\cap X_i\cap X_{j_k}$. 

\ \\ In total we have 

\begin{itemize}
    \item $c<f(\widetilde{y},s)$.

 \item $f(\widetilde{z},t)<c-\varepsilon_k$ for all $t\geq t_k$.

\item 
$\widetilde{y},\widetilde{z}\in \Omega_k\cap X_i$.

\end{itemize}

\noindent Now define:

 \begin{equation*}
   M_1:= \{x\in \Omega_k\cap X_i\,\,|\,\,for\,\,all\,\, t\in[t_k,t_{k+1}]\,\,we\,\,have\,\,f(x,t)\leq c-\varepsilon_{k+1} \,\},  
 \end{equation*}
\begin{equation*}
 M_2:= \{x\in \Omega_k\cap X_i\,\,|\,\,
  there\,\,exists\,\, t\in[t_k,t_{k+1}]\,\,such\,\,that\,\, 
c\leq f(x,t) \,\}.   
\end{equation*}

\ \\ \noindent By continuity of $f$ and by compactness of $[t_k,t_{k+1}]$, we have that $M_1$ and $M_2$ are closed, disjoints sets, which lie inside $\Omega_k\cap X_i$. Furthermore, $\widetilde{y}\in M_2$, and since $0<\varepsilon_{k+1}<\varepsilon_k$ we also have $\widetilde{z}\in M_1$. Hence, since $\Omega_k\cap X_i$ is connected, there exist $w\in (\Omega_k\cap X_i)\setminus(M_1\cup M_2)$. 

\ \\ That is, there exists $w\in \Omega_k\cap X_i$ such that 

\begin{equation*}\label{equation*14.11}
for\,\,all\,\,t\in[t_k,t_{k+1}]\,\,we\,\,have\,\,that\,\,f(w,t)<c,    
\end{equation*}
\begin{equation}\label{equation*15.11}
there\,\,exist\,\,s_{k+1}\in[t_k,t_{k+1}]\,\,with\,\,c-\varepsilon_{k+1}<f(w,s_{k+1}). 
\end{equation}\ \\ \noindent  By continuity of $f$ there exists $U_w$ a small enough neighborhood of $w$ such that both inequalities still hold in $U_w$. By connectedness via $(X_n)$, there exists $\Omega_{k+1}$ with compact closure such that 
$\overline{\Omega_{k+1}}\subset U_w\cap \Omega_k$ and $w\in\Omega_{k+1}$, and such that $\Omega_{k+1}$ is connected via $(X_n)$. So we have properties (A),(B),(C),(D).

\ \\ Define $j_{k+1}=i$. As $w\in X_{j_{k+1}}$, we get that $\Omega_{k+1}\cap X_{j_{k+1}}\neq\varnothing$. By the way we have chosen $t_{k+1}$, we have that for all $x\in\Omega_{k+1}\cap X_{j_{k+1}}$ and all $t\geq t_{k+1}$ we have that $f(x,t)<c-\varepsilon_{k+1}$, so we also get (E).

\ \\ As $t_k\leq s\leq t_{k+1}$ and $t_{k+1}\geq s+1$, we have that the sequences $(t_k)_{k=0}^{\infty}$ and $(s_k)_{k=1}^{\infty}$ are increasing and unbounded. This finishes the $k+1$ step.

\ \\ We now explain how to get uncountably many points $x_0$ satisfying $\limsup_{t\rightarrow\infty}f(x_0,t)=c$. 

\ \\ Assume by contradiction there exist at most countably many points $v$ satisfying that\\ $\limsup_{t\rightarrow\infty}f(v,t)=c$, and let $(v_k)_{k=0}^{\infty}$ be any enumeration of them. Now repeat the inductive construction until the point $w$ is chosen. Notice that as $w\in X_i$, by upper uniformity of $(X_n)$ we must have that $\limsup_{t\rightarrow\infty}f(w,t)=0$. So as $c>0$, we have that $w\neq v_k$ for all $k$. Now  choose $U_w$ to be a small enough neighborhood of $w$ such that both inequalities (\ref{equation*14.11}),(\ref{equation*15.11}) still hold in $U_w$ and such that $v_{k+1}\notin U_w$. So in particular, $v_{k+1}\notin \overline{\Omega_{k+1}}$.

\ \\ Repeating this argument at each step, we get that $v_k\notin\overline{\Omega_{k}}$ for all $k$. As the point $x_0$ we construct satisfies $x_0\in\cap_k \Omega_k$, we must have that $x_0\neq v_k$ for all $k$, contradicting the assumption that $(v_k)_{k=0}^{\infty}$ is an enumeration of all the points with $\limsup_{t\rightarrow\infty}f(v,t)=c$.\qed

\ \\ \indent Notice that the requirement that $f$ is continuous in both variables is essential in order to prove that the set $M_2$ is closed. Furthermore, recall that for $\gamma>0$ we have by Lemma \ref{Lemma6.111} that
$$
\limsup_{t\rightarrow\infty}\,\lambda_{\Theta,\psi_\gamma}(t)\,^{1+\gamma n/m}
=\limsup_{t\rightarrow\infty} \,\chi_{\gamma n/m}(\Theta,t),
$$

\noindent but that while $\lambda_{\Theta,\psi_\gamma}(t)$ is a continuous function of $\Theta$ and $t$, the function $\chi_{\gamma n/m}(\Theta,t)$ is only a piecewise continuous function in the $t$ variable.

\ \\ In the next chapter we prove Theorem 3 by applying Theorem 5 to the continuous function $\lambda_{\Theta,\psi}(t)$. A similar application to $\chi_{\gamma n/m}(\Theta,t)$ is not possible without some kind of a modification of Theorem 5 to a less restrictive version of it. 

\section{Applications; Proofs of Theorems 2 and 3}

\noindent We use Theorem 5 to prove Theorems 2 and 3. While both applications are very similar, we split between the Theorem 2 and Theorem 3 cases for extra clarity.

\subsection{Proof of Theorem 2}

Let $d\geq3$. Let $\{E_\varepsilon\,|\,\varepsilon\in (0,b]\,\}$ be a continuous decreasing exhaustion by compact sets of $\mathcal{X}_d$ for some $b\in(0,\infty]$. Let $m,n\in\mathbb{N}$ with $m+n=d$ and let $\Vec{\alpha}$ and $\Vec{\beta}$ be arbitrary weights on $\mathbb{R}^{d}$. 

\ \\ To prove Theorem 2, we apply Theorem 5 to the function $$f:M_{m,n}\times(0,\infty)\longrightarrow(0,b]$$
\begin{equation}\label{equation*24.444}
f(\Theta,t):=\inf\,\,\{\varepsilon>0\,\,|\,\, g_{t,\vec{\alpha},\vec{\beta}\,}\,\Lambda_{\Theta}\notin\,E_\varepsilon\,\}    
\end{equation} \noindent with the convention that $\inf\varnothing=b$.

\ \\ Notice that as $\{E_\varepsilon\,|\,\varepsilon\in(0,b]\,\}$ is a continuous decreasing exhaustion of $\mathcal{X}_d$ and that\\ $g_{t,\vec{\alpha},\vec{\beta}}$ \,acts continuously on $\mathcal{X}_d$, we have that $f(\Theta,t)$ is continuous.

\ \\ The application is as follows: First we show that there exists a sequence $(X_j)_{j=1}^{\infty}$ satisfying that $M_{m,n}$ is locally connected via $(X_j)$ and that $(X_j)$ satisfies the property of upper uniformity w.r.t. $f(\Theta,t)$ -- we do so while separating the cases $n=1$ and $n>1$. Then we explain why there exists a dense set satisfying $\{\Theta\in M_{m,n}\,|\,\limsup_{t\rightarrow\infty}f(\Theta,t)\geq b\,\}$. We do so by actually showing that $\{\Theta\in M_{m,n}\,|\,\limsup_{t\rightarrow\infty}f(\Theta,t)=b\,\}$ is a full Lebesgue measure set.  

\ \\ For the sequence $(X_j)$, we split between the cases of $n=1$ and $n\geq2$. Denote by $\norm{\cdot}_\infty$ the maximum norm.

\ \\ \indent 
 \underline{\pmb{Case 1, $n=1$.}} \,Notice that in this case the weight vector $\Vec{\beta}$ is just the number $1$. 

\ \\ For $\Vec{i}\in\mathbb{Z}^m$, $\Vec{z}\in\mathbb{Q}^m$ define $X_{\Vec{i},\Vec{z}}$ to be the line defined by $$X_{\Vec{i},\Vec{z}}:=\{y\,\vec{i}+\Vec{z}\,\,|\,\,y\in\mathbb{R}\},$$

\noindent and let $(X_j)_{j=1}^{\infty}$ be any enumeration of the set of lines $\{X_{\Vec{i},\Vec{z}}\,\,|\,\,\Vec{i}\in\mathbb{Z}^m,\Vec{z}\in\mathbb{Q}^m\,\}$.

\ \\ First we show $M_{m,1}$ is locally connected via $(X_j)$.

\ \\ Let $\vec{\Theta}\in M_{m,1}\cong\mathbb{R}^m$ and let $U_{\vec{\Theta}}$ be a neighborhood of $\vec{\Theta}$. We need to show that there exists $\Omega_{\vec{\Theta}}$ an open neighborhood of ${\vec{\Theta}}$ with compact closure which satisfies that $\overline{\Omega_{\vec{\Theta}}}\subset U_{\vec{\Theta}}$, there exists $j\in\mathbb{N}$ such that $\Omega_{\vec{\Theta}}\cap X_j\neq\varnothing$, and such that for any ${\vec{y}}\in \Omega_{\vec{\Theta}}$ and any $U_{\vec{y}}$ a neighborhood of ${\vec{y}}$, there exists $l\in\mathbb{N}$ such that $\Omega_{\vec{\Theta}}\,\cap\, X_l$ is connected, $\Omega_{\vec{\Theta}}\cap X_l\cap U_{\vec{y}}\neq\varnothing$ and $\Omega_{\vec{\Theta}}\cap X_l\cap X_j\neq\varnothing$.

\ \\ Choose $\Omega_{\vec{\Theta}}$ an open ball around $\Theta$ (w.r.t. the Euclidean norm, for example) such that $\overline{\Omega_{\vec{\Theta}}}\subset U_{\vec{\Theta}}$. By density of rational vectors, choose $j$ such that $\Omega_{\vec{\Theta}}\cap X_j\neq\varnothing$. Now let ${\vec{y}}\in\Omega_{\vec{\Theta}}$ and $U_{\vec{y}}$ a neighborhood of ${\vec{y}}$. Choose ${\vec{z}_0}\in\mathbb{Q}^m$ with ${\vec{z}_0}\in U_{\vec{y}}\cap\Omega_{\vec{\Theta}}$. So for any ${\vec{i}}\in\mathbb{Z}^m$ we have that ${\vec{z}_0}\in \Omega_{\vec{\Theta}}\cap X_{{\vec{i}},{\vec{z}_0}}\cap U_{\vec{y}}$. Let $\vec{z}_1\in\mathbb{Q}^m$ with $\vec{z}_1\in\Omega_{\vec{\Theta}}\cap X_j$ and with $\vec{z}_1\neq\vec{z}_0$. Such a vector exists because $\mathbb{Q}^m$ is dense in $X_j$. Now let $\vec{i}_0\in\mathbb{Z}^m$ be a multiple of $\vec{z}_1-\vec{z}_0$ so that $\vec{z}_1\in X_{\vec{i}_0,\vec{z}_0}$. Then the line $X_{\vec{i}_0,\vec{z}_0}$ intersects $\Omega_{\vec{\Theta}}$ in a segment, which is a connected set, and this segment intersects both $\Omega_{\vec{\Theta}}\cap X_j$ and $\Omega_{\vec{\Theta}}\cap U_{\vec{y}}$. So for $X_l=X_{{\vec{i}}_0,{\vec{z}_0}}$ we get that $M_{m,1}$ is locally connected via $(X_j)$.

\ \\ In order to show $(X_j)$ satisfies the property of upper uniformity, we use the following claim, which holds for both $n=1$ and $n\geq2$:

\begin{Claim}
Let $(X_j)\subset M_{m,n}$ such that for any $c>0$, any $j$ and any $\Theta\in X_j$ there exists $T_0$ such that for any $t\geq T_0$ there exists $\vec{v}_t\in\Lambda_\Theta$ with $\lVert \,g_{t,\vec{\alpha},\vec{\beta}}\,\,\Vec{v}_{t}\lVert_{\infty}\,<c$. Then $(X_j)$ satisfies the property of upper uniformity w.r.t. $f(\Theta,t)$ for $f$ as in (\ref{equation*24.444}). 
\end{Claim}

\noindent \textbf{Proof of Claim 1.} \,\,\,To show that $(X_j)$ satisfies the property of upper uniformity w.r.t. $f(\Theta,t)$, it is enough to show that for any $\eta>0$, any $j$ and any $\Theta\in X_j$ there exists $T_0$ such that for any $t\geq T_0$ we have that $ g_{t,\vec{\alpha},\vec{\beta}\,}\,\Lambda_{\Theta}\notin\,E_\eta$. Notice that this is a stronger property than upper uniformity, in which
the above needs to hold only within a compact subset of $X_j$.

\ \\ Recall that by Mahler compactness criterion and Dirichlet's theorem, $\{K_{\varepsilon}(\norm{\cdot}_{\infty})\,|\,\varepsilon\in[0,1]\,\}$ 
is a continuous decreasing exhaustion of $\mathcal{X}_d$ by compact sets. That is, there exists $c\in(0,1)$ such that $E_\eta\subset K_{c}(\norm{\cdot}_{\infty})$. So its enough to show that there exists $T_0$ such that for any $t\geq T_0$ we have that $g_{t,\vec{\alpha},\vec{\beta}}\,\Lambda_\Theta\notin K_{c}(\norm{\cdot}_{\infty})$. So by the definition of $K_{c}(\norm{\cdot}_{\infty})$, we get the claim. \qed

\ \\ Now denote $\alpha_{max}:=\norm{\Vec{\alpha}}_{\infty}$ and  $i_{max}:=\lVert\Vec{i\,}\lVert_{\infty}$. 

\ \\ Let $\Vec{x}=y\,\vec{i}+\Vec{z}\in X_{\Vec{i},\Vec{z}}$\, for some $\Vec{i},\Vec{z}$ and $y\in\mathbb{R}$, and let $c>0$. Let $z_0\in\mathbb{N}$ be any natural number such that $z_0\vec{z}\in\mathbb{Z}^m$.

\ \\  By Dirichlet's theorem, let $(q_k,p_k)_{k=0}^{\infty}\in\mathbb{N}\times\mathbb{Z}$ be any sequence with $1=q_0<q_1<q_2<\,...$ and $|q_0y-p_0|>|q_1y-p_1|>|q_2y-p_2|>\,...$ \,such that $q_{k+1}|q_ky-p_k|\leq1$. 

\ \\ Define also $q_{-1}:=q_0=1$ and $p_{-1}:=p_0$. So we have $q_{k+1}|q_ky-p_k|\leq1$ for every $k\geq-1$. 

\ \\ Let $(t_k)_{k=0}^{\infty}$ be the increasing sequence which satisfies for all $k\geq-1$ that
$$t_{k+1}^{-1}\,q_{k+1}=i_{max}\,t_{k+1}^{\,\alpha_{max}}\,|q_ky-p_k\,|$$ and define $t_{-1}:=0$. 

\ \\ For $t>0$ such that $t\in(t_k,t_{k+1}]$ define

\begin{equation*}
f^*_y(t):=\max\, \Biggl\{ 
\begin{aligned}
&\,\,\frac{q_k}{t}
\\
&\,\, i_{max}\,t^{\,\alpha_{max}}\,|q_ky-p_k\,|.
\end{aligned} 
\end{equation*}

\ \\ Furthermore, for $t\in(t_k,t_{k+1}]$  define  $\vec{l}_{t}:=z_0(q_k\Vec{z}+p_k\Vec{i}\,)\in\mathbb{Z}^m$, and define $$\Vec{v}_{t}:=\begin{pmatrix}
I_m & \Vec{x} \\
0 & 1 \\
\end{pmatrix}\begin{pmatrix}
\Vec{l}_{t}  \\
-z_0q_k \\
\end{pmatrix}=-z_0\big((q_ky-p_k)\vec{i}\,,\,q_k\,\big)\in\mathbb{R}^{m+1}$$

\noindent noticing that $\Vec{v}_{t}\in\Lambda_{\Vec{x}}$ for all $t>0$.

\ \\ So for all $t>0$ we have
\begin{equation*}\label{equation*13}\begin{split} 
g_{t,\vec{\alpha},1}\,\Vec{v}_{t}=-z_0\Big( t^{\alpha_1}i_1(q_ky-p_k)\,,\,...\,,\,t^{\alpha_m}i_m(q_ky-p_k)\,,\,q_kt^{-1}\Big)
\end{split}
\end{equation*}

\noindent and for all $t>0$ we have that  

\begin{equation*}\label{equation*28.22}
\lVert g_{t,\vec{\alpha},1}\,\Vec{v}_{t}\lVert_{\infty}\,\,\leq z_0\,f^*_y(t).
\end{equation*}

\ \\ So by Claim 1 we get that in order to show that $(X_j)$ has the upper uniformity w.r.t. $f(\Theta,t)$ in case $n=1$, it is enough to show that for any $c>0$ there exists $T_0>0$ such that $f^*_y(t)<z_0^{-1}c$ for any $t>T_0$, with $T_0$ independent of $y$.

\ \\ Note that $f^*_y$ obtains all its local maxima at $(t_k)_{k=0}^{\infty}$. 

Within these times we have that
\begin{equation*}
f^*_y(t_{k+1})^2=i_{max}\,t_{k+1}^{\,\alpha_{max}-1}    \,q_{k+1}|q_ky-p_k|\leq i_{max}\,t_{k+1}^{\,\alpha_{max}-1}.     
\end{equation*}

\noindent

\ \\ So by the definition of $f^*_y$ we have that for all $t$ with $t>t_0$ we have that $$f^*_y(t)\leq \,\sqrt{i_{max}}\,\,t^{\,(\alpha_{max}-1)/2},$$ where as $m\geq2$ we have that $\alpha_{\max}<1$. So choose $T_1$ such that $ z_0^{-1}\,\sqrt{i_{max}}\,T_1^{\,(\alpha_{max}-1)/2}<c$. Note that although the sequence $(t_k)_{k=0}^{\infty}$ depends on $y$, the fixed time $T_1$ does not. Also note that $t_0$ is the first local maxima of $f^*_y$. 

\ \\ If $t_0\leq T_1$ then we set $T_0:=T_1$ and we are done. 

\ \\ If $T_1<t_0$ then first notice that by definition of $T_1$ we have that $f^*_y(t)<c$ for all $t\geq t_0$. \ \\ In this case, we also have:

\begin{itemize}
    \item If $f^*_y$ is increasing in $(0,t_0)$ then $f^*_y(T_1)< f^*_y(t)<f^*_y(t_0)<c$ for all $t\in(T_1,t_0)$. So again we set $T_0:=T_1$ and we are done. 
    
    \item  If $f^*_y$ is decreasing in $(0,s)$ for some $s\in(0,t_0)$ then notice that $f^*_y(t)=1/t$ for all $t\in(0,s]$. In this case, we also have: 
    \begin{itemize}
        \item If $1/c< s$ then $f^*_y(t)<c$ for all $t\in(1/c,s]$. As $f^*_y$ is increasing in $(s,t_0)$ and as we assume $T_1<t_0$, we also have $f^*_y(t)<f^*_y(t_0)<c$ for all $t\in(s,t_0)$. So we set $T_0=1/c$ and we are done.

        \item If $s\leq 1/c$, then first assume by contradiction that $1/c<t_0$. Then $c\leq f^*_y(1/c)\leq f^*_y(t_0)$, contradicting our assumption that $f^*_y(t_0)<c$. So we must have $t_0\leq 1/c$, and so we set again $T_0=1/c$ \,and we are done.  
         
    \end{itemize}
\end{itemize}
\begin{equation}\label{list2}
\end{equation}

\noindent So we get in total that by setting $T_0=\max\,\{T_1,1/c\}$ we have that $(X_j)$ has the upper uniformity w.r.t $f(\Theta,t)$ in case $n=1$.

\ \\ \indent \underline{\pmb{Case 2, $n\geq2$.}} \,Case 2 is easier. For $\Vec{i}\in\mathbb{Q}^n$ and $\Vec{z}\in\mathbb{Q}^m$ define $Y_{\Vec{i},\Vec{z}}$ to be the $m(n-1)$ dimensional affine plane defined by

$$Y_{\Vec{i},\Vec{z}}:=\{\Theta\in M_{m,n}\,|\,\Theta\Vec{i}=\Vec{z}\,\},$$

\noindent and let $(Y_j)_{j=1}^{\infty}$ be any enumeration of the set $\{ Y_{\Vec{i},\Vec{z}}\,|\,\vec{i}\in\mathbb{Q}^n,\,\vec{z}\in\mathbb{Q}^m\,\}$. Note that $\cup_{j}Y_j$ contains all rational matrices.

\ \\ To show 
$M_{m,n}$ is connected via $(Y_j)$, let $\Theta\in M_{m,n}$ and let $U_{\Theta}$ be a neighborhood of $\Theta$.

\ \\ Choose again $\Omega_{\Theta}$ an open ball around $\Theta$ (w.r.t. the Euclidean norm on $M_{m,n}$ identified as $\mathbb{R}^{mn}$, for example) such that $\overline{\Omega_\Theta}\subset U_\Theta$. Notice that for any $l\in\mathbb{N}$ we have that $\Omega_\Theta\cap Y_l$ is connected. By density of rational matrices, choose $j$ such that $\Omega_\Theta\cap Y_j\neq\O$. Now let $y\in\Omega_\Theta$ and $U_y$ a neighborhood of $y$. Choose $Q_0\in M_{m,n}(\mathbb{Q})$ with $Q_0\in U_y\cap\,\Omega_\Theta$. Notice that for all $\vec{i}_0\in\mathbb{Q}^n$ there exists $\vec{z}_0\in\mathbb{Q}^m$ with $Q_0\in \Omega_{\Theta}\cap Y_{{\vec{i}_0},{\vec{z}_0}}\cap U_y$. Choose ${\vec{i}}_0,\vec{z}_0$ such that $\Omega_\Theta\cap Y_{{\vec{i}}_0,{\vec{z}_0}}\cap Y_j\neq\varnothing$. Such a choice exists as a solution of a nonhomogeneous system of two   linear equations and as we assume $n\geq2$. So for $Y_l=Y_{{\vec{i}}_0,{\vec{z}_0}}$ we get that $M_{m,n}$ is locally connected via $(Y_j)$.

\ \\ For upper uniformity, let $\Theta\in Y_{\Vec{i},\Vec{z}}$\, for some $\Vec{i},\Vec{z}$. Let $z_0\in\mathbb{N}$ such that $z_0\vec{i}\in\mathbb{Z}^n$, $z_0\vec{z}\in\mathbb{Z}^m$. So the vector  $$\Vec{v}_{t}:=\begin{pmatrix}
I_m & \Theta \\
0 & I_n \\
\end{pmatrix}\begin{pmatrix}
-z_0\Vec{z}\,  \\
z_0\Vec{i}\,
\end{pmatrix}=z_0(\Theta\Vec{i}-\Vec{z}\,,\,\Vec{i}\,)=z_0(\Vec{0}\,,\,\Vec{i}\,)\in\mathbb{R}^{m+n}$$ satisfies that $\Vec{v}_{t}\in\Lambda_{\Theta}$, and for all $t>0$ we have

\begin{equation*}g_{t,\vec{\alpha},\vec{\beta}}\,\,\Vec{v}_{t}=z_0\,(\,\Vec{0}\,,\,t^{-\beta_1}i_1,\,...\,,t^{-\beta_n}i_n\,) .  
\end{equation*}

\ \\ So again by Claim 1 we get that $(X_j)$ has the upper uniformity w.r.t. $f(\Theta,t)$ in case $n>1$. 

\ \\ \ \\ \indent We now explain why there exists a full Lebesgue measure set with $\limsup_{t\rightarrow\infty}f(\Theta,t)=b$ for any $\Theta$ in this set. 

\ \\ That is, we need to show that for any $c<b$ there is an unbounded  positive sequence $(t_k)$ such that $g_{t_k,\vec{\alpha},\vec{\beta}}\,\Lambda_\Theta\in E_c$ for Lebesgue almost every $\Theta$.  

\ \\ We do so by using the same argument as in \cite{KR3}, \cite{KW} (Theorem 1.4 in both),  using the following equidistribution theorem of Kleinbock and Weiss (\cite{KW}, Theorem 2.2):

\begin{Theorem}
Let $\vec{\alpha},\vec{\beta}$ arbitrary weights on $\mathbb{R}^d$. Let $F\in C_c(\mathcal{X}_d), B\subset M_{m,n}$ be bounded with positive Lebesgue measure,
and $\delta>0$ be given. Then there exists $t_0>0$ such that for all $t>t_0$,

$$ |\,vol(B)^{-1}\int_B F(g_{t,\vec{\alpha},\vec{\beta}}\,\Lambda_\Theta)\,dvol(\Theta)-\int_{\mathcal{X}_d}F(L)\,dm_{\mathcal{X}_d}(L)\,|<\delta. $$

\noindent Here, the integrals are taken with respect to the Lebesgue measure $vol$ on $M_{m,n}$ and the Haar-Siegel probability measure $m_{\mathcal{X}_d}$ on $\mathcal{X}_d$.
\end{Theorem} 

\noindent Now for $i\in\mathbb{N}$ define $$B_i:=\bigcap_{\,t>i}\,\{\Theta\in M_{m,n}\,|\,g_{t,\vec{\alpha},\vec{\beta}}\,\Lambda_\Theta\notin E_c\},$$ \noindent and assume by contradiction that $vol(B_i)>0$. Choose $B\subset B_i$ compact with positive measure as well. Take a non-negative function $F\in C_c(\mathcal{X}_d)$ which is supported on $E_c$ but vanishes outside of it, and define $\delta:=\frac{1}{2}\int_{\mathcal{X}_d}F\,dm_{\mathcal{X}_d}>0$. Applying Theorem 6 with $t>i$ sufficiently large we get that $\int_B F\,(g_{t,\vec{\alpha},\vec{\beta}}\Lambda_\Theta)\,dvol(\Theta)>0$, contradicting the definition of $B_i$. 

\ \\ So $vol\,(\bigcup_{i\in\mathbb{N}}B_i)=0$, i.e. for Lebesgue almost every $\Theta$ there is an unbounded  positive sequence $(t_k)$ such that $g_{t_k,\vec{\alpha},\vec{\beta}}\,\Lambda_\Theta\in E_c$.

 \ \\ This finishes the proof of Theorem 2. \qed

\subsection{Proof of Theorem 3}

Let $m$ and $n$ such that $\max\,(m,n)>1$, let $\norm{\cdot}_m$ and $\norm{\cdot}_n$ be arbitrary norms on $\mathbb{R}^m$ and $\mathbb{R}^n$, and let $\Delta$ be the minimal constant which satisfies (\ref{equation*4}). Let $\psi(t)$ be a positive continuous decreasing function such that $\psi(t)=o(t^{-1})$, where in case $n=1$ we assume additionally that $\psi(t)^{-1/m}=o(t)$.

\ \\ To prove Theorem 3, we apply Theorem 5 to the function
$$g:M_{m,n}\times(0,\infty)\longrightarrow(0,b]$$
\begin{equation}\label{equation*31.333}
g(\Theta,t):=\lambda_{\Theta,\psi}(t)    
\end{equation} 

\noindent for $\lambda_{\Theta,\psi}(t)$ as in ($\ref{equation*5.1111})$.

\ \\ Notice that as $\psi(t)$ is continuous and by discreteness of $\mathbb{Z}^d$ we have that $g(\Theta,t)$ is continuous.

\ \\ The application is the same as in the proof of Theorem 2: First we show that there exists a sequence $(X_j)_{j=1}^{\infty}$ satisfying that $M_{m,n}$ is connected via $(X_j)$ and that $(X_j)$ has the upper uniformity w.r.t. $g(\Theta,t)$ -- we do so while separating the cases $n=1$ and $n>1$, and by using the same sequences $(X_j)$ for both cases as in the proof of Theorem 2. Then we explain why there exists a dense set satisfying $\{\Theta\in M_{m,n}\,|\,\limsup_{t\rightarrow\infty}g(\Theta,t)= \infty\,\}$.

\ \\ Denote by $\widetilde{\mathcal{X}_d}:=GL_d(\mathbb{R})/GL_d(\mathbb{Z})$ the space of lattices in $\mathbb{R}^d$ (not necessarily unimodular).

\ \\ Given $t>0$ define
$$\pmb{a_{\,\psi,t}}:=diag\,(\psi(t)^{-1/m},\,...\,,\psi(t)^{-1/m},t^{-1/n},\,...\,,t^{-1/n})\in GL_d(\mathbb{R}).$$ 

\ \\ For the sequence $(X_j)$, we split again between the cases of $n=1$ and $n\geq2$.

\ \\  \indent 
 \underline{\pmb{Case 1, $n=1$.}} \,Notice that in this case the weight vector $\Vec{\beta}$ is just the number $1$, for $\Vec{i}\in\mathbb{Z}^m$, 
$\vec{z}\in\mathbb{Q}^m$ define as before $X_{\Vec{i},\Vec{z}}$ to be the line defined by $$X_{\Vec{i},\Vec{z}}:=\{y\,\vec{i}+\Vec{z}\,\,|\,\,y\in\mathbb{R}\},$$

\noindent and let $(X_j)_{j=1}^{\infty}$ be any enumeration of the set of lines $\{X_{\Vec{i},\Vec{z}}\,\,|\,\,\Vec{i},\Vec{z}\in\mathbb{Z}^m\,\}$. 

\ \\ By section 5.1, $M_{m,1}$ is connected via $(X_j)$.

\ \\ In order to show $(X_j)$ has the upper uniformity, we use the following claim, which again holds for both $n=1$ and $n\geq2$:

\begin{Claim}
Let $(X_j)\subset M_{m,n}$ such that for any $c>0$, any $j$ and any $\Theta\in X_j$ there exists $t_0$ such that for any $t\geq t_0$ there exists $\vec{v}_t\in\Lambda_\Theta$ with $\lVert a_{\psi,t}\,\Vec{v}_{t}\lVert_{\infty}\,<c$. Then $(X_j)$ has the upper uniformity w.r.t. $g(\Theta,t)$ for $g$ as in (\ref{equation*31.333}).
\end{Claim}

\noindent \textbf{Proof of Claim 2.} \,\,\,To show that $(X_j)$ has the upper uniformity w.r.t. $g(\Theta,t)$, it is enough to show that for any $\eta>0$, any $j$ and any $\Theta\in X_j$ there exists $t_0$ such that for any $t\geq t_0$ we have that $ \lambda_{\Theta,\psi}(t)<\eta$.  Notice that this is a stronger property than upper uniformity, in which
the above needs to hold only within a compact subset of $X_j$. 

\ \\ By definition of $\lambda_{\Theta,\psi}(t)$ we have that $\lambda_{\Theta,\psi}(t)\leq \norm{a_{\psi,t}\,\vec{v}_t}_{m,n}$ for any $\vec{v}_t\in\Lambda_\Theta$. So by norm equivalence, we get the claim. \qed

\ \\ Denote as before  $i_{max}:=\lVert\Vec{i\,}\lVert_{\infty}$. 

\ \\ Let $\Vec{x}=y\,\vec{i}+\Vec{z}\in X_{\Vec{i},\Vec{z}}$\, for some $\Vec{i},\Vec{z}$ and $y\in\mathbb{R}$, and let $c>0$. Let $z_0\in\mathbb{N}$ be any natural number such that $z_0\vec{z}\in\mathbb{Z}^m$.

\ \\  By Dirichlet's theorem, let $(q_k,p_k)_{k=0}^{\infty}\in\mathbb{N}\times\mathbb{Z}$ be any sequence with $1=q_0<q_1<q_2<\,...$ and $|q_0y-p_0|>|q_1y-p_1|>|q_2y-p_2|>\,...$ \,such that $q_{k+1}|q_ky-p_k|\leq1$. 

\ \\ Define also $q_{-1}:=q_0=1$ and $p_{-1}:=p_0$. So we have $q_{k+1}|q_ky-p_k|\leq1$ for every $k\geq-1$. 

\ \\ Let $(t_k)_{k=0}^{\infty}$ be the increasing sequence which satisfies for all $k\geq-1$ that
$$t_{k+1}^{-1}\,q_{k+1}=i_{max}\,\psi(t_{k+1})^{\,-1/m}\,|q_ky-p_k\,|$$ and define $t_{-1}:=0$. 

\ \\ For $t>0$ such that $t\in(t_k,t_{k+1}]$ define

\begin{equation*}
g^*_y(t):=\max\, \Biggl\{ 
\begin{aligned}
&\,\,\frac{q_k}{t}
\\
&\,\, i_{max}\,\psi(t)^{\,-1/m}\,|q_ky-p_k\,|.
\end{aligned} 
\end{equation*}

\ \\ Furthermore, for $t\in(t_k,t_{k+1}]$  define  $\vec{l}_{t}:=z_0(q_k\Vec{z}+p_k\Vec{i}\,)\in\mathbb{Z}^m$, and define $$\Vec{v}_{t}:=\begin{pmatrix}
I_m & \Vec{x} \\
0 & 1 \\
\end{pmatrix}\begin{pmatrix}
\Vec{l}_{t}  \\
-z_0q_k \\
\end{pmatrix}=-z_0\big((q_ky-p_k)\vec{i}\,,\,q_k\,\big)\in\mathbb{R}^{m+1}$$

\noindent noticing that $\Vec{v}_{t}\in\Lambda_{\Vec{x}}$ for all $t>0$.

\ \\ So for all $t>0$ we have
\begin{equation*}\label{equation*132.1}\begin{split} 
a_{\psi,t}\,\Vec{v}_{t}=-z_0\Big(\psi(t)^{-1/m}i_1(q_ky-p_k)\,,\,...\,,\,\psi(t)^{-1/m}i_m(q_ky-p_k)\,,\,q_kt^{-1}\Big).
\end{split}
\end{equation*}

\noindent and for all $t>0$ we have that  

\begin{equation*}\label{equation*28.22.1}
\lVert a_{\psi,t}\,\Vec{v}_{t}\lVert_{\infty}\,\,\leq z_0\,g^*_y(t).
\end{equation*}

\noindent So by Claim 2 we get that in order to show that $(X_j)$ has the upper uniformity w.r.t. $g(\Theta,t)$ in case $n=1$, it is enough to show that for any $c>0$ there exists $T_0>0$ such that $g^*_y(t)<z_0^{-1}c$ for any $t>T_0$, with $T_0$ independent of $y$.

\ \\ Note that $g^*_y$ obtains all its local maxima at $(t_k)_{k=0}^{\infty}$. 

Within these times we have that
\begin{equation*}
g^*(y,t_{k+1})^2=i_{max}\,\dfrac{\psi(t_{k+1})^{-1/m}}{t_{k+1}}\,q_{k+1}|q_ky-p_k|\leq i_{max}\,\dfrac{\psi(t_{k+1})^{-1/m}}{t_{k+1}}.      
\end{equation*}
\noindent

\ \\ So by the definition of $g^*_y$ we have that for all $t$ with $t>t_0$ we have that 

$$g^*_y(t)\leq \,\sqrt{i_{max}}\,\Big(\dfrac{\psi(t_{k+1})^{-1/m}}{t_{k+1}}\Big)^{1/2},$$ 

\noindent and we assume  $\psi(t)^{-1/m}=o(t)$. So choose $T_1$ such that $ z_0^{-1}\,\sqrt{i_{max}}\,\big(\psi(T_1)^{-1/m}\,T_1^{-1}\big)^{1/2}<c$. Note that although the sequence $(t_k)_{k=0}^{\infty}$ depends on $y$, the fixed time $T_1$ does not. Also note that $t_0$ is the first local maxima of $g^*_y$. 

\ \\ Now define $T_0=\max\,\{T_1,1/c\}$. By the same reasoning as in (\ref{list2}) we get that $(X_j)$ has the upper uniformity w.r.t $g(\Theta,t)$ in case $n=1$.

\ \\ \indent \underline{\pmb{Case 2, $n\geq2$.}} \,For $\Vec{i}\in\mathbb{Q}^n$ and $\Vec{z}\in\mathbb{Q}^m$ define $Y_{\Vec{i},\Vec{z}}$ to be the $m(n-1)$ dimensional affine plane defined by

$$Y_{\Vec{i},\Vec{z}}:=\{\Theta\in M_{m,n}\,|\,\Theta\Vec{i}=\Vec{z}\,\},$$

\noindent and let $(Y_j)_{j=1}^{\infty}$ be any enumeration of the set $\{ Y_{\Vec{i},\Vec{z}}\,|\,\vec{i}\in\mathbb{Z}^n,\,\vec{z}\in\mathbb{Z}^m\,\}$. Note that $\cup_{j}Y_j$ contains all rational matrices.

\ \\ Again, $M_{m,n}$ is connected via $(Y_j)$ by the proof given in section 5.1.

\ \\ For upper uniformity, let $\Theta\in Y_{\Vec{i},\Vec{z}}$\, for some $\Vec{i},\Vec{z}$. Let $z_0\in\mathbb{N}$ such that $z_0\vec{i}\in\mathbb{Z}^n$, $z_0\vec{z}\in\mathbb{Z}^m$. So the vector  $$\Vec{v}_{t}:=\begin{pmatrix}
I_m & \Theta \\
0 & I_n \\
\end{pmatrix}\begin{pmatrix}
-z_0\Vec{z}\,  \\
z_0\Vec{i}\,
\end{pmatrix}=z_0(\Theta\Vec{i}-\Vec{z}\,,\,\Vec{i}\,)=z_0(\Vec{0}\,,\,\Vec{i}\,)\in\mathbb{R}^{m+n}$$ satisfies that $\Vec{v}_{t}\in\Lambda_{\Theta}$, and for all $t>0$ we have

\begin{equation*}a_{\psi,t}\,\,\Vec{v}_{t}=z_0\,(\,\Vec{0}\,,\,t^{-1/n}i_1,\,...\,,t^{-1/n}i_n\,) .  
\end{equation*}

\ \\ So again by Claim 2 we get that $(X_j)$ has the upper uniformity w.r.t. $g(\Theta,t)$ in case $n>1$. 

\ \\ \ \\ \indent We now explain why there exists a full Lebesgue measure set with $\limsup_{t\rightarrow\infty}g(\Theta,t)=\infty$ for any $\Theta$ in this set. 

\ \\ We use Theorem 6 as in section 5.1 on the decreasing exhaustion  $\{K_\varepsilon(\norm{\cdot}_{m,n})\,|\,\varepsilon\in(0,r_{\norm{\cdot}_{m,n}}]\,\}$ for $r_{\norm{\cdot}_{m,n}}>0$ as in (\ref{equation*10.1111}) and get that for Lebesgue almost every $\Theta\in M_{m,n}$ and every $c$ with $ c<r_{\norm{\cdot}_{m,n}}$ there is an unbounded  positive sequence $(t_k)$ such that $g_{t_k,\vec{m}_1,\vec{n}_1}\,\Lambda_\Theta\in K_c(\norm{\cdot}_{m,n})$.

\ \\ So by (\ref{equation*10.111}) and (\ref{equation*16.33}) we have that $\limsup_{t\rightarrow\infty}\lambda_{\Theta,t^{-1}}(t)=r_{\norm{\cdot}_{m,n}}$ for Lebesgue almost every $\Theta\in M_{m,n}$, and as we assume $\psi(t)=o(t^{-1})$ we get by Lemma \ref{Lemma2} that $\limsup_{t\rightarrow\infty}g(\Theta,t)=\infty$ for Lebesgue almost every $\Theta\in M_{m,n}$. 

\ \\ This finishes the proof of Theorem 3. \qed

\ \\ \textbf{Remark.} Recalling that $\mathbb{D}_{1,1}\neq[0,\Delta]=[0,1]$, a natural question for ending is why do we have such a different phenomenon of the Dirichlet spectrum for $n=m=1$ from a Theorem 5 perspective.  

\ \\ To see why our arguments break down for $(m,n)=(1,1)$, note that in this case the sets $(X_j)$ we used are just points and the condition of $\mathbb{R}$ being locally connected via $(X_j)$ does not hold.

\section{Technical lemmas}

\noindent We start by stating Lemmas 1-4. Notice that while Lemmas 1 and 2 are used in previous sections of the paper, Lemmas 3 and 4 are used only in order to prove Lemmas 1 and 2.
We use Lemma 1 during the proof of Theorem 3 in order to prove that $\limsup_{t\rightarrow\infty}g(\Theta,t)=\infty$ for Lebesgue almost every $\Theta\in M_{m,n}$, for $g$ as in (\ref{equation*31.333}). We use Lemma 2 in order to deduce Corollary 1 from Theorem 3.

\ \\ Let $\psi(t)$ be a positive continuous decreasing function with $\lim_{t\rightarrow\infty}\psi(t)=0$, $\Theta\in M_{m,n}$ fixed, and let $\norm{\cdot}_m$ and $\norm{\cdot}_n$ be arbitrary norms on $\mathbb{R}^m$ and $\mathbb{R}^n$. By discreteness of $\mathbb{Z}^d$ and since $\lambda_{\Theta,\psi}$ is continuous and $\psi$ is decreasing to zero, there exists a sequence $(t_k,\vec{q}_k,\vec{p}_k)\subset\mathbb{R}_{>0}\times\mathbb{Z}^n\times\mathbb{Z}^m$ with $1=t_0<t_1<t_2...$ such that if $t\in[t_{k},t_{k+1}]$  then 

\begin{equation}\label{equation*33.111}
\lambda_{\Theta,\psi}(t)=\max\, \Biggl\{
\begin{array}{l}
t^{-1/n}\,\norm{\Vec{q}_k}_n
   \vspace*{0.2cm}\\
\psi(t)^{-1/m}\,\norm{\Theta\Vec{q}_k-\vec{p}_k\,}_m.
\end{array}    
\end{equation}

\ \\ In case there exists $(\vec{Q},\vec{P}\,)\in\mathbb{Z}^n\times\mathbb{Z}^m$ such that $\Theta\vec{Q}=\vec{P}$ then the sequence $(t_k,\vec{q}_k,\vec{p}_k)$ is finite, and for all $t$ sufficiently large we have that $\lambda_{\Theta,\psi}(t)=t^{-1/n}\,\lVert \vec{Q}\lVert_n$. Otherwise, the sequence $(t_k,\vec{q}_k,\vec{p}_k)$ is infinite, and we have that $t_k\longrightarrow\infty$. 

\ \\ In both of these cases, we say $\boldsymbol{(t_k,\Vec{q}_k,\Vec{p}_k)}$ \textbf{realizes} $\boldsymbol{\lambda_{\Theta,\psi}}$.

\ \\ We now state the lemmas proved in this section:

\begin{Lemma}\label{Lemma2}
Let $\norm{\cdot}_m$ and $\norm{\cdot}_n$ be arbitrary norms on $\mathbb{R}^m$ and $\mathbb{R}^n$, and let $\psi(t)$ be a positive continuous decreasing function such that $\psi(t)=o(t^{-1})$. Let $\Theta\in M_{m,n}$ such that $\limsup_{t\rightarrow\infty}\,\lambda_{\Theta,t^{-1}}(\Theta,t)>0$. Then $$\limsup_{t\rightarrow\infty}\,\lambda_{\Theta,\psi}(\Theta,t)=\infty.$$
\end{Lemma}

\begin{Lemma}\label{Lemma6.111}
Let $\norm{\cdot}_m$ and $\norm{\cdot}_n$ be arbitrary norms on $\mathbb{R}^m$ and $\mathbb{R}^n$, and let $\gamma>0$. Then $$\limsup_{t\rightarrow\infty}\,\lambda_{\Theta,\psi_\gamma}(t)\,^{1+\gamma n/m}
=\limsup_{t\rightarrow\infty} \,\chi_{\gamma n/m}(\Theta,t).$$
\end{Lemma}

\begin{Lemma}\label{Lemma3.1} Let $\norm{\cdot}_m$ and $\norm{\cdot}_n$ be arbitrary norms on $\mathbb{R}^m$ and $\mathbb{R}^n$, let $\psi(t)$ be a positive 
continuous decreasing function with $\lim_{t\rightarrow\infty}\psi(t)=0$, and let $\Theta\in M_{m,n}$. Let $(t_k,\vec{q}_k,\vec{p}_k)\subset\mathbb{R}_{>0}\times\mathbb{Z}^n\times\mathbb{Z}^m$ which realizes $\lambda_{\Theta,\psi}$ as in (\ref{equation*33.111}). Then
$$ 0<\norm{\Vec{q}_0}_n<\norm{\Vec{q}_1}_n<...<\norm{\Vec{q}_k}_n$$
$$\norm{\Theta\Vec{q}_0-\Vec{p}_0}_m>\norm{\Theta\Vec{q}_1-\Vec{p}_1}_m>\,...\,>\norm{\Theta\Vec{q}_k-\Vec{p}_k}_m.
$$     

\noindent In case the sequence is infinite, we also have $\lim\limits_{k\rightarrow\infty}\norm{\Vec{q}_k}_n=\infty$ and $\lim\limits_{k\rightarrow\infty}\norm{\Theta\Vec{q}_k-\Vec{p}_k}_m=0$.
\end{Lemma}

\begin{Lemma}\label{Lemma3}
Let $\norm{\cdot}_m$ and $\norm{\cdot}_n$ be arbitrary norms on $\mathbb{R}^m$ and $\mathbb{R}^n$, let $\psi(t)$ be a positive 
continuous decreasing function with $\lim_{t\rightarrow\infty}\psi(t)=0$, and let $\Theta\in M_{m,n}$. Let $(t_k,\vec{q}_k,\vec{p}_k)\subset\mathbb{R}_{>0}\times\mathbb{Z}^n\times\mathbb{Z}^m$ which realizes $\lambda_{\Theta,\psi}$ as in (\ref{equation*33.111}). 

\ \\ Then there exists a sequence of times $(z_k)_{k=1}^{\infty}$ with $z_k\in(t_k,t_{k+1})$ such that for all $k$ we have that $\lambda_{\Theta,\psi}(t)$ is decreasing in $[t_k,z_k]$, increasing in $[z_k,t_{k+1}]$,
and we have that 

\begin{itemize}
\item $t\in[t_k,z_k]\Longrightarrow \lambda_{\Theta,\psi}(t)=\lVert t^{-1/n}\,\Vec{q}_k\lVert_n$
\item $t\in[z_k,t_{k+1}]\Longrightarrow \lambda_{\Theta,\psi}(t)=\norm{\psi(t)^{-1/m}\,(\Theta\Vec{q}_k-\Vec{p}_k)}_m$.
\end{itemize}
\end{Lemma}

\subsection{Proofs}

\ \\ \textbf{Proof of Lemma 3.} \,\,\, Let $k\in\mathbb{N}$. We can assume that the segment in which $(\Vec{q}_k,\Vec{p}_k)$ realizes $\lambda_{\Theta,\psi}$ is of maximum of length. I.e. if $t\notin[t_k,t_{k+1}]$ then $$\lambda_{\Theta,\psi}(t)<\max\, \Biggl\{
\begin{array}{l}
t^{-1/n}\,\norm{\Vec{q}_k}_n
   \vspace*{0.2cm}\\
\psi(t)^{-1/m}\,\norm{\Theta\Vec{q}_k-\vec{p}_k\,}_m
\end{array}. $$

\ \\ Let $S,T$ with $t_k<S<t_{k+1}<T<t_{k+2}$. Then we have:

\begin{align*}
&\lambda_{\Theta,\psi}(S)=\max \Biggl\{
\begin{array}{l}
S^{-1/n}\,\lVert\Vec{q}_k\lVert_n
   \vspace*{0.2cm}\\
\psi(S)^{-1/m}\,\lVert\Theta\Vec{q}_k-\vec{p}_k\,\lVert_m
\end{array}
&&<
\,\,\, \max \Biggl\{
\begin{array}{l}
S^{-1/n}\,\lVert\Vec{q}_{k+1}\lVert_n
   \vspace*{0.2cm}\\
\psi(S)^{-1/m}\,\lVert\Theta\Vec{q}_{k+1}-\vec{p}_{k+1}\,\lVert_m
\end{array} 
\\[10pt]
&\lambda_{\Theta,\psi}(t_{k+1})=\max \Biggl\{
\begin{array}{l}
t_{k+1}^{-1/n}\,\lVert\Vec{q}_k\lVert_n
   \vspace*{0.2cm}\\
\psi(t_{k+1})^{-1/m}\,\lVert\Theta\Vec{q}_k-\vec{p}_k\,\lVert_m
\end{array}
&&=\,\,\, \max \Biggl\{
\begin{array}{l}
t_{k+1}^{-1/n}\,\lVert\Vec{q}_{k+1}\lVert_n
   \vspace*{0.2cm}\\
\psi(t_{k+1})^{-1/m}\,\lVert\Theta\Vec{q}_{k+1}-\vec{p}_{k+1}\,\lVert_m
\end{array} 
\\[10pt]
&\lambda_{\Theta,\psi}(T)=\max \Biggl\{
\begin{array}{l}
T^{-1/n}\,\lVert\Vec{q}_{k+1}\lVert_n
   \vspace*{0.2cm}\\
\psi(T)^{-1/m}\,\lVert\Theta\Vec{q}_{k+1}-\vec{p}_{k+1}\,\lVert_m
\end{array}
&&<\,\,\,\max \Biggl\{
\begin{array}{l}
T^{-1/n}\,\lVert\Vec{q}_k\lVert_n
   \vspace*{0.2cm}\\
\psi(T)^{-1/m}\,\lVert\Theta\Vec{q}_k-\vec{p}_k\,\lVert_m
\end{array}.
\end{align*}

\ \\ As we assume that $\psi$ is decreasing to zero, we get Lemma 3. \qed

\ \\   

\noindent \textbf{Proof of Lemma 4.} \,\,\,  By Lemma \ref{Lemma3.1} we have for all $k$ 

\begin{equation*}\label{equation*35.11}
\lambda_{\Theta,\psi}(t_{k+1})=t_{k+1}^{-1/n}\,\lVert\Vec{q}_{k+1}\lVert_n\,>\,\psi(t_{k+1})^{-1/m}\,\lVert\Theta\Vec{q}_{k+1}-\vec{p}_{k+1}\,\lVert_m
\end{equation*}

\noindent and so $\lambda_{\Theta,\psi}$ is decreasing in some right neighbourhood of $t_{k+1}$.

\ \\ As $\psi$ is decreasing to zero, the function $h(t):= \max \Biggl\{
\begin{array}{l}
t^{-1/n}\,\lVert\Vec{q}_{k+1}\lVert_n
   \vspace*{0.2cm}\\
\psi(t)^{-1/m}\,\lVert\Theta\Vec{q}_{k+1}-\vec{p}_{k+1}\,\lVert_m
\end{array}$ 

\ \\ is the maximum between two positive functions -- the first is decreasing to zero, the second increases to infinity -- hence there exists a minimum point $z\in(t_{k+1},\infty)$ such that $h(t)$ is decreasing in $(0,z]$ and increases to infinity in $[z,\infty)$.  

\ \\ So we are only left with showing that $z<t_{k+2}$. Indeed, assume by contradiction $z\geq t_{k+2}$. Then we have

 \begin{equation}\label{equation*38.1111}
\lambda_{\Theta,\psi}(t_{k+2})=t_{k+2}^{-1/n}\,\lVert\Vec{q}_{k+1}\lVert_n\,\geq\,\psi(t_{k+2})^{-1/m}\,\lVert\Theta\Vec{q}_{k+1}-\vec{p}_{k+1}\,\lVert_m.
\end{equation}

\noindent So as $(t_k,\Vec{q}_k,\Vec{p}_k)$ realizes $\lambda_{\Theta,\psi}$, we must have that 

\begin{equation*}
t_{k+2}^{-1/n}\,\lVert\Vec{q}_{k+1}\lVert_n\,=  \max \Biggl\{
\begin{array}{l}
t_{k+2}^{-1/n}\,\lVert\Vec{q}_{k+2}\lVert_n
   \vspace*{0.2cm}\\
\psi(t_{k+2})^{-1/m}\,\lVert\Theta\Vec{q}_{k+2}-\vec{p}_{k+2}\,\lVert_m.
\end{array}   
\end{equation*}

\noindent By Lemma \ref{Lemma3.1} $,\lVert\Vec{q}_{k+1}\lVert_n<\lVert\Vec{q}_{k+2}\lVert_n$ so we must have 

$$t_{k+2}^{-1/n}\,\lVert\Vec{q}_{k+1}\lVert_n\,=\psi(t_{k+2})^{-1/m}\,\lVert\Theta\Vec{q}_{k+2}-\vec{p}_{k+2}\,\lVert_m,$$

\noindent hence by (\ref{equation*38.1111})
$$ \psi(t_{k+2})^{-1/m}\,\lVert\Theta\Vec{q}_{k+2}-\vec{p}_{k+2}\,\lVert_m\,\geq\,\psi(t_{k+2})^{-1/m}\,\lVert\Theta\Vec{q}_{k+1}-\vec{p}_{k+1}\,\lVert_m,$$

\ \\ which cannot hold by Lemma \ref{Lemma3.1}. \qed 

\ \\ 

\noindent \textbf{Proof of Lemma 1.} \,\,\, Notice that as we assume that $\limsup_{t\rightarrow\infty}\,\lambda_{\Theta,t^{-1}}(\Theta,t)>0$ we have that the sequence which realizes $\lambda_{\Theta,t^{-1}}$ is infinite, hence also the sequence which realizes $\lambda_{\Theta,\psi}$. Now assume by contradiction that $\limsup_{t\rightarrow\infty}\lambda_{\Theta,\psi}(t)< \infty$, \,and let $(t_k,\vec{q}_k,\vec{p}_k)_{k=0}^{\infty}\subset\mathbb{R}_{>0}\times\mathbb{Z}^n\times\mathbb{Z}^m$ which realizes $\lambda_{\Theta,\psi}$. 

\ \\ By Lemma \ref{Lemma3}, $(t_k)_{k=1}^{\infty}$ is the sequence of local maxima of $\lambda_{\Theta,\psi}$ and for all $k\geq0$ we have that
\begin{equation}\label{equation*25.11}\begin{split}
\lambda_{\Theta,\psi}(t_{k+1})=\lVert\psi(t_{k+1})^{-1/m}\,(\Theta\Vec{q}_k-\Vec{p}_k)\lVert_{m
}
=\lVert t_{k+1}^{-1/n}\,\Vec{q}_{k+1}\lVert_{n}.  
\end{split}
\end{equation}

\noindent So by (\ref{equation*25.11}) we have that \begin{equation*}
\begin{split}
\infty>  \limsup_{t\rightarrow\infty}\lambda_{\Theta,\psi}(t)^{m+n}&=\limsup_{k\rightarrow\infty}\lambda_{\Theta,\psi}(t_{k+1})^{m+n}
\\&=\limsup_{k\rightarrow\infty}\,(t_{k+1}\,\psi(t_{k+1}))^{-1}\,
\lVert \Vec{q}_{k+1}\lVert_n^{n}\,\,\lVert(\Theta\Vec{q}_k-\Vec{p}_k)\lVert_m^{m}.
\end{split}
\end{equation*}

\noindent As we assume $\psi(t)=o(t^{-1})$ we must have

\begin{equation}\label{equation*31.1}
\limsup_{k\rightarrow\infty}\,\lVert \Vec{q}_{k+1}\lVert_n^{n}\,\,\lVert(\Theta\Vec{q}_k-\Vec{p}_k)\lVert_m^{m}\,=0. 
\end{equation}

\ \\ Now let $(s_k)_{k=1}^{\infty}$ be the sequence which satisfies for all $k\geq0$ that 

\begin{equation*}
\lVert s_{k+1}^{-1/n}\,\Vec{q}_{k+1}\lVert_n=\lVert s_{k+1}^{1/m}\,(\Theta\Vec{q}_k-\Vec{p}_k)\lVert_m.
\end{equation*}

\noindent By Lemma \ref{Lemma3.1} we have $s_1<s_2<...<s_k<...$ and $s_k\longrightarrow\infty$.

\ \\ Now define a function $g:[s_1,\infty)\longrightarrow(0,\infty)$ such that for $t\in[s_k,s_{k+1})$ we set

\begin{equation*}
g(t):=\max\Biggl\{
\begin{array}{l}
t^{-1/n}\,\lVert\Vec{q}_k\lVert_n
   \vspace*{0.2cm}\\
t^{1/m}\,\lVert\Theta\Vec{q}_k-\vec{p}_k\,\lVert_m
\end{array}
\end{equation*}

\noindent Notice that there exists a sequence of times $(z_k)_{k=1}^{\infty}$ with $z_k\in(s_k,s_{k+1})$ such that for all $k\geq1$ we have that $g$ is decreasing in $[s_k,z_k]$, increasing in $[z_k,s_{k+1}]$,
and we have that 

\begin{itemize}
\item $t\in[s_k,z_k]\Longrightarrow g(t)=\lVert t^{-1/n}\,\Vec{q}_k\lVert_n$
\item $t\in[z_k,s_{k+1}]\Longrightarrow g(t)=\lVert t^{1/m}\,(\Theta\Vec{q}_k-\Vec{p}_k)\lVert_m$.
\end{itemize}

\noindent Indeed, assume by contradiction that $g$ is increasing in some right neighborhood of $s_{k+1}$ for some $k$. So we must have that $\lVert s_{k+1}^{-1/n}\,\Vec{q}_{k+1}\lVert_n\leq\lVert s_{k+1}^{1/m}\,(\Theta\Vec{q}_{k+1}-\Vec{p}_{k+1})\lVert_m$, which cannot hold since $\lVert\Theta\Vec{q}_k-\Vec{p}_k\lVert_m>\lVert\Theta\Vec{q}_{k+1}-\Vec{p}_{k+1}\lVert_m$ and the way we chose the sequence $(s_k)_{k=1}^{\infty}$.

\ \\ Similarly, assume by contradiction that $g$ is decreasing in some left neighborhood of $s_{k+1}$ for some $k$. So we must have that $\lVert s_{k+1}^{1/m}\,(\Theta\Vec{q}_{k}-\Vec{p}_{k})\lVert_m\leq\lVert s_{k+1}^{-1/n}\,\Vec{q}_{k}\lVert_n$, which cannot hold since $\lVert \Vec{q}_k\lVert_n<\lVert \Vec{q}_{k+1}\lVert_n$ and the way we chose the sequence $(s_k)_{k=1}^{\infty}$.

\ \\ So $(s_k)_{k=1}^{\infty}$ is the sequence of local maxima of $g(t)$, and for all $t\geq s_1$ we have $\lambda_{\Theta,t{-1}}\leq g(t)$. 

\ \\ In particular, \begin{equation*}
\begin{split}
\limsup_{t\rightarrow\infty}\lambda_{\Theta,t^{-1}}(t)^{m+n}&\leq\limsup_{t\rightarrow\infty}g(t)^{m+n}=
\limsup_{k\rightarrow\infty}g(s_{k+1})^{m
+n}
\\&=\limsup_{k\rightarrow\infty}\,\,
\big(\lVert s_{k+1}^{-1/n}\,\Vec{q}_{k+1}\lVert_n\big)^n\,\big(\lVert s_{k+1}^{1/m}\,(\Theta\Vec{q}_k-\Vec{p}_k)\lVert_m\big)^m
\\ &=\limsup_{k\rightarrow\infty}\,
\lVert \,\Vec{q}_{k+1}\lVert_n^{n}\,\,\lVert(\Theta\Vec{q}_k-\Vec{p}_k)\lVert_m^{m}\,=0
\end{split}
\end{equation*}

\noindent where the last equality follows from (\ref{equation*31.1}).

\ \\ So we get in total $\limsup_{t\rightarrow\infty}\lambda_{\Theta,t^{-1}}(t)=0$, contradicting our main assumption on $\Theta$.

\ \\ This finishes the proof of Lemma 1. \qed

\ \\ 

\noindent \textbf{Proof of Lemma 2.} \,\,\,In case $\Theta\vec{q}\in\mathbb{Z}^m$ \,for some $\vec{q}\in\mathbb{Z}^n\setminus\{\vec{0}\}$ we have that both quantities are equal to zero (in this case, $\Theta\in X_j$ for some $j\in\mathbb{N}$ and $(X_j)$ as in the proof of Theorem 2 in section 5, where we explain why $\limsup_{t\rightarrow\infty}\,\lambda_{\Theta,\psi_\gamma}(t)=0$ for any $j\in\mathbb{N}$ and any $\Theta\in X_j$). So assume $\Theta\vec{q}\notin\mathbb{Z}^m$ for all non-zero $\vec{q}\in\mathbb{Z}^n$.

\ \\ Assume that $\limsup_{t\rightarrow\infty}\lambda_{\Theta,\psi_\gamma}(t)^{1+\gamma n/m}= \,L$.  

\ \\ Let $(t_k,\vec{q}_k,\vec{p}_k)_{k=0}^{\infty}\subset\mathbb{R}_{>0}\times\mathbb{Z}^n\times\mathbb{Z}^m$ which realizes $\lambda_{\Theta,\psi_\gamma}$ as in (\ref{equation*33.111}). So if $t\in[t_{k},t_{k+1}]$ we have that \begin{equation*}
\lambda_{\Theta,\psi_\gamma}(t)=\max \Biggl\{
\begin{array}{l}
t^{-1/n}\,\norm{\Vec{q}_k}_n
   \vspace*{0.2cm}\\
t^{\,\gamma/m}\,\norm{\Theta\Vec{q}_k-\vec{p}_k\,}_m\end{array}.    
\end{equation*}

\ \\ Again by Lemmas \ref{Lemma3.1} and \ref{Lemma3} we have that $(t_k)_{k=1}^{\infty}$ is the sequence of local maxima of $\lambda_{\Theta,\psi_\gamma}$ and for all $k\geq0$ we have that \begin{equation*}\label{equation*40.11}
\begin{split}
\lambda_{\Theta,\psi_\gamma}(t_{k+1})=\lVert t_{k+1}^{\,\,\gamma/m}\,(\Theta\Vec{q}_k-\Vec{p}_k)\lVert_{m
}=\lVert t_{k+1}^{-1/n}\,\Vec{q}_{k+1}\lVert_{n}.  
\end{split}
\end{equation*}
   
\noindent So in particular we have that 
\begin{equation*}\label{equation*41.11}
\begin{split}
t_{k+1}=\big(\,\lVert \,\Vec{q}_{k+1}\lVert_{n} \,\,\lVert\,(\Theta\Vec{q}_k-\Vec{p}_k)\lVert_{m
}^{-1}\big)^{nm/(\gamma n+m)},
\end{split}
\end{equation*}

\noindent and we have that \begin{equation}\label{equation*59.111}
\begin{split}
\limsup_{t\rightarrow\infty}\lambda_{\Theta,\psi_\gamma}(t)^{2}&=\limsup_{k\rightarrow\infty}\lambda_{\Theta,\psi_\gamma}(t_{k+1})^{2}
\\&=\limsup_{k\rightarrow\infty}\,t_{k+1}^{^{\,\gamma/m-1/n}}\,
\lVert \Vec{q}_{k+1}\lVert_n\,\,\lVert(\Theta\Vec{q}_k-\Vec{p}_k)\lVert_m
\\&=\limsup_{k\rightarrow\infty}\,t_{k+1}^{^{\,(\gamma n-m)/nm}}\,
\lVert \Vec{q}_{k+1}\lVert_n\,\,\lVert(\Theta\Vec{q}_k-\Vec{p}_k)\lVert_m
\\&=\limsup_{k\rightarrow\infty}\,\big(\,\lVert \,\Vec{q}_{k+1}\lVert_{n}\,\lVert(\Theta\Vec{q}_k-\Vec{p}_k)\lVert_{m
}^{-1}\big)^{(\gamma n-m)/(\gamma n+m)}\lVert \Vec{q}_{k+1}\lVert_n\,\,\lVert(\Theta\Vec{q}_k-\Vec{p}_k)\lVert_m
\\&=\limsup_{k\rightarrow\infty}\,\lVert \Vec{q}_{k+1}\lVert_{n}^{\,2\gamma n/(m+\gamma n)}\,\lVert(\Theta\Vec{q}_k-\Vec{p}_k)\lVert_{m
}^{\,2m/(m+\gamma n)}
\end{split}
\end{equation}

\noindent and as we assume $\limsup_{t\rightarrow\infty}\lambda_{\Theta,\psi_\gamma}(t)^{1+\gamma n/m}= \,L$ we have in total that

\begin{equation*}\label{equation*60.111}
\limsup_{k\rightarrow\infty}\,\lVert \Vec{q}_{k+1}\lVert_{n}^{\gamma n/m}\,\lVert(\Theta\Vec{q}_k-\Vec{p}_k)\lVert_{m
}=L.
\end{equation*}

\ \\ Now define a function $f:[\norm{\vec{q}_0},\infty)\longrightarrow[0,\infty)$ by

\begin{equation*}
f(t):=t^{\gamma n/m}\,\lVert\Theta\vec{q}_k-\vec{p}\,\lVert_m\,\,\,\,\,\,\,\
\end{equation*} \noindent for $t\in[\,\norm{\vec{q}_k},\norm{\vec{q}_{k+1}}\,)$. 

\ \\ So in particular we have
\begin{equation*}
\begin{split}
\limsup_{t\rightarrow\infty}\chi_{\gamma n/m}(t)\leq\limsup_{t\rightarrow\infty}f(t)=\limsup_{k\rightarrow\infty}\,\lVert \vec{q}_{k+1}\lVert_n^{\gamma n/m}\,\lVert\Theta\vec{q}_k-\vec{p}\,\lVert_m
=L.
\end{split}
\end{equation*}

\ \\ Now assume by contradiction that $\limsup_{t\rightarrow\infty}\chi_{\gamma n/m}(t)<L$, \,and let $(\vec{Q}_k,\vec{P}_k)_{k=0}^{\infty}\subset\mathbb{Z}^n\times\mathbb{Z}^m$ which realizes $\chi_{\gamma n/m}(\Theta,t)$. I.e. if $t\in[\lVert\vec{Q}_{k}\lVert_n,\lVert\vec{Q}_{k+1}\lVert_n)$ then we have that \begin{equation*}
\chi_{\gamma n/m}(\Theta,t)=t^{\gamma n/m}\,\lVert\Theta\vec{Q}_k-\vec{P}_k\lVert_m, 
\end{equation*}

\noindent and by our assumption we have that

\begin{equation}\label{equation*64}
\limsup\limits_{t\rightarrow\infty}\chi_{\gamma n/m}(\Theta,t)=\limsup\limits_{k\rightarrow\infty}\,\lVert\vec{Q}_{k+1}\lVert_n^{\gamma n/m}\,\,\lVert\Theta\vec{Q}_k-\vec{P}_k\lVert_m<L.
\end{equation}

\noindent Now repeat the process as in the proof of Lemma \ref{Lemma2}.

\ \\ Let $(s_k)_{k=1}^{\infty}$ be the sequence which satisfies that for all $k\geq0$ we have 

\begin{equation*}
\lVert s_{k+1}^{-1/n}\,\Vec{Q}_{k+1}\lVert_n=\lVert s_{k+1}^{\gamma/m}\,(\Theta\Vec{Q}_k-\Vec{P}_k)\lVert_m.
\end{equation*}

\noindent So as we have that $\Theta\vec{q}\notin\mathbb{Z}^m$ for all non-zero $\vec{q}$, we have that $s_1<s_2<...<s_k<...$ and that $s_k\longrightarrow\infty$.

\ \\ Now define a function $g:[s_1,\infty)\longrightarrow(0,\infty)$ such that for $t\in[s_k,s_{k+1})$ we set

\begin{equation*}
g(t):=\max\Biggl\{
\begin{array}{l}
t^{-1/n}\,\lVert\Vec{Q}_k\lVert_n
   \vspace*{0.2cm}\\
t^{\gamma/m}\,\lVert\Theta\Vec{Q}_k-\vec{P}_k\,\lVert_m
\end{array}
\end{equation*}

\ \\ By the same arguing as in Lemma \ref{Lemma2}, $(s_k)_{k=1}^{\infty}$ is the sequence of local maxima of $g(t)$, and for all $t\geq s_1$ we have that $\lambda_{\Theta,\psi_\gamma}\leq g(t)$. 

\ \\ In particular, by the same calculations as in (\ref{equation*59.111}) we have that

\begin{equation*}
\begin{split}
L^{2m/(m+\gamma n)}&=\limsup_{t\rightarrow\infty}\lambda_{\Theta,\psi_\gamma}(t)^{2}\leq\limsup_{t\rightarrow\infty}g(t)^{2}=\limsup_{k\rightarrow\infty}g(s_{k+1})^2\\&=\limsup_{k\rightarrow\infty}\,\lVert \Vec{Q}_{k+1}\lVert_{n}^{(2\gamma n)/(m+\gamma n)}\,\lVert(\Theta\Vec{Q}_k-\Vec{P}_k)\lVert_{m
}^{(2m)/(m+\gamma n)}.
\end{split}
\end{equation*}

\noindent So in total we have

\begin{equation*}
L^m=\limsup_{k\rightarrow\infty}\,\lVert \Vec{Q}_{k+1}\lVert_{n}^{\gamma n}\,\lVert(\Theta\Vec{Q}_k-\Vec{P}_k)\lVert_{m
}^{m}
\end{equation*}

\noindent contradicting formula (\ref{equation*64}).

\ \\ This finishes the proof of Lemma 2. \qed

\newpage
\bibliographystyle{alpha}
\bibliography{ref}

\end{document}